# A THEORY OF BOND PORTFOLIOS

By Ivar Ekeland and Erik Taflin

*University of British Columbia and EISTI*

We introduce a bond portfolio management theory based on foundations similar to those of stock portfolio management. A general continuous-time zero-coupon market is considered. The problem of optimal portfolios of zero-coupon bonds is solved for general utility functions, under a condition of no-arbitrage in the zero-coupon market. A mutual fund theorem is proved, in the case of deterministic volatilities. Explicit expressions are given for the optimal solutions for several utility functions.

**1. Introduction.** This paper is a first step toward a unified theory of portfolio management, including both stocks and bonds. There is a gap between the traditional approaches to manage bond portfolios and stock portfolios. Managing bond portfolios relies on concepts such as duration, sensibility and convexity, while managing stock portfolios relies on optimization of expected utility. We give two results toward bridging this gap. First, we set up and solve the problem of managing a bond portfolio by optimizing (over all self-financing trading strategies for a given initial capital) the expected utility of the final wealth. Second, we express the solution of this problem as portfolios of self-financing trading strategies which include naturally stocks and bonds.

The well-established theory of portfolio management, initiated in the seminal papers [13, 14, 15], [20] and further developed by many, see [12], [17] and references therein, does not apply as it stands to bond portfolios. The difficulty here is that stocks and bonds differ in many ways, the most important of which is the fact that bonds mature at a prescribed date (time of maturity) after which they disappear from the market, whereas the characteristics of a stock do not change, except in reaction to business news or management decisions. Another difference is that in an unconstrained market, the time of maturity can take an infinity of values, so there is an infinity of different bonds. As a first consequence, the price of a stock depends only on the risks it carries (market risk, idiosyncratic risk), whereas the price of a bond depends both on the risks it carries (interest rate risk, credit risk) and on time to maturity. Mathematically, this is expressed by the fact that the stochastic differential equations used to model stock prices are usually autonomous (meaning that the coefficients are time-independent functions of the prices, as in geometric Brownian motion or mean-reverting processes), whereas any model for bond prices must incorporate the fact that the volatility goes to zero when time to maturity goes to zero. So the mathematical analysis of a portfolio including stocks and bonds is complicated by the fact that the prices for each type of assets evolve according to different rules, even in the most elementary case. An added difficulty, due to the maturity dependence, is that certain strategies which are possible for stocks are no longer allowable for bonds: a simple buy-and-hold strategy, for instance, results in converting bonds to cash on maturity. The particular case of strategies involving only a finite number of bonds, all with maturities exceeding the portfolio management horizon, is similar to the case of a pure stock market (with stochastic interest rate). Optimal portfolios for such cases were obtained in [10]. An optimal portfolio problem in a truly maturity-dependent context of reinsurance contracts was solved in [21] for discrete time.

Our suggestion is to work in a "moving frame," that is, to consider time to maturity, instead of maturity, as the basic variable on which the zero-coupon depends at each time. At time $t$, there will be a curve $S \to p_t(S)$, $S \geq 0$, where $p_t(S)$ is the price of a standard zero-coupon maturing at time $t + S$. Here $S$ is time to maturity and $T = t + S$ is time of maturity. Such a parameterization was introduced in [16]. When $t$ changes, so does the curve $p_t$, and a bond portfolio then is simply a linear functional operating on the space of such curves. Now from the financial point of view, this can be seen in different perspectives: (1) The static point of view, say, is to consider the portfolio at time $t$ simply as a linear combination (possibly infinite) of standard zero-coupons, each of which has a fixed time of maturity $T \geq t$. Such a portfolio has to be rebalanced each time a zero-coupon in the portfolio comes to maturity. (2) The dynamic point of view is to consider the portfolio at time $t$ as a linear combination of self-financing instruments, each one with a fixed time to maturity $S \geq 0$. We term such an instrument a Roll-Over and it is simply a certain $t$-dependent multiple of a zero-coupon with time to maturity $S$, independent of $t$ (see Remark 2.7). Its price has a simple expression, given by (2.33). Such instruments were introduced earlier in [19] under the name "rolling-horizon bond." Roll-Overs behave like stocks, in the sense that their time to maturity is constant through time, so that their price depends only on the risk they carry. One can then envision a program where portfolios are expressed as combinations of stocks and Roll-Overs, which are treated in a uniform fashion.

However, it is well known that this program entails mathematical difficulties. The first one is that rewriting the equations for bond prices in the moving frame introduces the operator $\frac{\partial}{\partial S}$, which has to make sense as an unbounded operator in the space $H$ of curves $p_t$ chosen to describe zero-coupon prices. The second one is that this space $H$ has to be contained in the space of all continuous functions on $\mathbb{R}^+$, so that its dual $H^*$ contains the Dirac masses $\delta_{T-t}$, corresponding at time $t$ to one zero-coupon of maturity $T$, but should not be too small, otherwise $H^*$ will contain many more objects which cannot easily be interpreted as bond portfolios. In this paper we choose $H$ to be a standard Sobolev space, which in particular is a Hilbert space. Bond portfolios are then simply elements of the Hilbert space $H^*$. Reference [1] introduced portfolios being signed finite Borel measures. They also are elements of $H^*$. The analysis is in our case simplified by the fact that $H$ and $H^*$ are Hilbert spaces. In a different context, Hilbert spaces of forward rates were considered in [3] and [6]. The image of these spaces, under the nonlinear map of forward rates to zero-coupons prices, is locally included in $H$.

We believe that this abstract, Hilbertian approach opens up many possibilities. In this paper, as mentioned above we explore one, namely portfolio management. We give existence theorems for very general utility functions and for $H$-valued price processes driven by a cylindrical Wiener process, that is, in our case by a countable number of independent Brownian motions. We give explicit solutions, taking advantage of the Hilbertian setting. These solutions are expressed in terms of (nonunique) combinations of classical zero-coupon bonds [i.e., financial interpretation (1) above], but the optimal strategy can readily be translated in terms of Roll-Overs, which may not be marketed, although they are self-financing [i.e., financial interpretation (2) above]. If the price of bonds depends on a $d$-dimensional Brownian motion, then the optimal strategy can be expressed as a linear combination of $d$ bonds and in certain cases these can be any $d$ marketed bonds, with time of maturity exceeding the time horizon of the optimal portfolio problem.

The outline of the paper is as follows. We begin by setting up the appropriate framework in Section 2, where bond portfolios are defined as elements of a certain Hilbert space $H^*$. Bond dynamics are prescribed in (2.11) according to the HJM methodology [7] and a self-financing portfolio is defined (cf. [1]) by formulas (2.27) and (2.28). An arbitrage-free market is prescribed according to Condition A and we introduce certain self-financing trading strategies with fixed time to maturity, which we call Roll-Overs (Remark 2.7). The optimal portfolio problem is set up in Section 3, and solved in two special cases, the first being when the underlying Brownian motion is finite-dimensional (Theorem 3.6), the second being when it is infinite-dimensional, but the market price of risk is a deterministic function of time (Theorem 3.8). Examples of closed-form solutions are then given in Section 4. All our portfolios are functions of the market price of risk, similar to those giving the Merton portfolio in the case of stocks. This indicates that our treatment indeed unifies bond and stock portfolio management.

Mathematical proofs are provided in Section 5 and the Appendix. In the Appendix we state and prove some existence results and estimates for infinite-dimensional processes with stochastic volatility that we have not found in standard references such as [4] or [9]. We note that the appropriate mathematical framework for the study of infinite-dimensional (cylindrical) processes is the theory of Hilbert–Schmidt operators, to which we appeal in the proofs, although we have avoided it in the statement of the results.

Several remarks of a mathematical nature are made in Section 5. Remark 5.1 justifies our market condition (Condition A), in Remark 5.4 it is shown that our results apply to certain incomplete markets and in Remark 5.5 a Hamilton–Jacobi–Bellman approach is considered. We note, in Remark 5.6, that our existence result for certain utility functions with asymptotic elasticity equal to 1 stands in apparent contrast with the earlier result of [11] and [12] for stock portfolios. This is because we have used a narrower definition (Condition A) of arbitrage-free prices for bonds.

## 2. The bond market.
We consider a continuous-time bond market and without restriction we can assume that only zero-coupon bonds are available. The time horizon in our model is some finite date $\bar{T} > 0$. At any date $t \in \mathbb{T} = [0, \bar{T}]$, one can trade zero-coupon bonds with maturity $s \in [t, \infty[$. Bonds with maturity $s = t$ at time $t$ will be assimilated to money in a current account (see (ii) of Example 2.6 and cf. [2]).

Uncertainty is modeled by a filtered probability space $(\Omega, P, \mathcal{F}, \mathcal{A})$; here $\mathcal{A} = \{\mathcal{F}_t | 0 \leq t \leq \bar{T}\}$ is a filtration of the $\sigma$-algebra $\mathcal{F}$. The random sources are given by independent Brownian motions $W^i$, $i \in \mathbb{I}$. The index set $\mathbb{I}$ can be finite, $\mathbb{I} = \{1, \ldots, \bar{m}\}$, or infinite, $\mathbb{I} = \mathbb{N}^* = \mathbb{N} - \{0\}$. The filtration $\mathcal{A}$ is generated by the $W^i$, $i \in \mathbb{I}$.

2.1. *Zero-coupons and state space.* As usual, we denote by $B(t, s)$ the price at time $t$ of a zero-coupon bond yielding one unit of account at time $s$, $0 \leq t < s$, so that $B(t, t) = 1$. It is an $\mathcal{F}_t$-measurable random variable. In order to introduce interest rates let us assume that, almost surely, the function $s \mapsto B(t, s)$ is strictly positive and $C^1$. We denote by $r(t)$ the spot interest rate at $t$ and by $f_t(S)$ the instantaneous forward rate contracted at $t \in \mathbb{T}$ for time to maturity $S$:

$$(2.1) \qquad r(t) = f_t(0) \quad \text{and} \quad f_t(S) = -\frac{1}{B(t, t+S)} \frac{\partial B}{\partial S}(t, t+S),$$

which is allowed to be negative. $\bar{B}(t,s)$ denotes the price discounted to time 0:

$$(2.2) \qquad \bar{B}(t,s) = B(t,s) \exp\left(-\int_0^t r(\tau)\,d\tau\right).$$

It will be convenient to characterize zero-coupon bonds by their time to maturity. For this reason we introduce the $\mathcal{A}$-adapted $C^1([0,\infty[)$-valued processes $p$ and $\bar{p}$ defined by

$$(2.3) \qquad p_t(S) = B(t,t+S) \quad \text{and} \quad \bar{p}_t(S) = \bar{B}(t,t+S),$$

where $t \in \mathbb{T}$ and $S \geq 0$. This parameterization was introduced in [16]. One should here take care that $S$ is the time to maturity and not the maturity itself. Note that $p_t(0) = 1$. We shall call $p_t$ (resp. $\bar{p}_t$) the zero-coupon bond (resp. discounted zero-coupon bond) state at time $t$. For simplicity we will also use zero-coupon bond state or just state for both cases. The state at time $t$ can thus be thought of as the curve: zero-coupon bonds price at the instant $t$ as function of time to maturity. Obviously

$$(2.4) \qquad B(t,s) = p_t(s-t) \quad \text{and} \quad \bar{B}(t,s) = \bar{p}_t(s-t),$$

where $t \in \mathbb{T}$ and $s - t \geq 0$.

More generally we will assume the processes $p$ and $\bar{p}$ to take values in a certain Sobolev space $H$, the zero-coupon bond state space. Our choice of $H$ is motivated by the following considerations:

(a) $H$ is a space of continuous functions going to zero at infinity, because zero-coupon bond prices are continuous with respect to time to maturity and they tend to zero as time to maturity tends to infinity.

(b) $H$ should be a Hilbert space, because it is the simplest possible infinite-dimensional topological vector space.

Conditions (a) and (b) leave us little choice, except to take $H$ to be a Sobolev space such as $H^s(]0,\infty[)$, with $s > 1/2$ (see below). Note that further conditions should be required for the model to be completely realistic:

(c) $p_t(S)$ must be differentiable with respect to $S$ at $S = 0$, so that the spot interest rate is well defined.

(d) $p_t(S)$ should be positive for all $S > 0$ and $p_t(0) = 1$.

(e) $p_t(S)$ should be decreasing with respect to $S$.

Conditions (c) and (d) will be satisfied as a result of our model. However, to include simple Gaussian interest rate models, we will not impose condition (e). The state space of portfolios at each time $H^*$, which is the dual of the zero-coupon bond state space, will contain measures as it shall. If wanted, we can now choose $H$ such that portfolios have certain regularity properties, for example, such that derivatives of measures are not elements of $H^*$. We next define $H$ and recall certain elementary facts concerning Sobolev spaces.

For $s \in \mathbb{R}$, let $H^s = H^s(\mathbb{R})$ (cf. Section 7.9 of [8]) be the usual Sobolev space of real tempered distributions $f$ on $\mathbb{R}$ such that the function $x \mapsto (1 + |x|^2)^{s/2} \hat{f}(x)$ is an element of $L^2(\mathbb{R})$, where $\hat{f}$ is the Fourier transform [in $\mathbb{R}^n$ we denote $x \cdot y = \sum_{1 \le i \le n} x_i y_i$, $x, y \in \mathbb{R}^n$, and we define the Fourier transform $\hat{f}$ of $f$ by $\hat{f}(y) = (2\pi)^{-n/2} \int_{\mathbb{R}^n} \exp(-iy \cdot x) f(x) \, dx$] of $f$, endowed with the norm:

$$\|f\|_{H^s} = \left( \int (1 + |x|^2)^s |\hat{f}(x)|^2 \, dx \right)^{1/2}.$$

All the $H^s$ are Hilbert spaces. Clearly, $H^0 = L^2$ and $H^s \subset H^{s'}$ for $s \ge s'$ and in particular $H^s \subset L^2 \subset H^{-s}$, for $s \ge 0$. If $f$ is $C^n$, $n \in \mathbb{N}$, and if $f$ together with its $n$ first derivatives belong to $L^2$, then $f \in H^n$. For every $s$, the space $C_0^\infty(\mathbb{R})$ of $C^\infty$ functions with compact support is dense in $H^s$. For every $s > 1/2$, by the Sobolev embedding theorems, we have $H^s \subset C^0 \cap L^\infty$. In addition $H^s$ is a Banach algebra for $s > 1/2$: if $f \in H^s$ and $g \in H^s$, then $fg \in H^s$ and the multiplication is continuous. Also, if $s > 1/2$, $f \in H^s$ and $g \in H^{-s}$, then $fg \in H^{-s}$ and the multiplication is continuous also here.

We define, for $s \in \mathbb{R}$, a continuous bilinear form on $H^{-s} \times H^s$ by

$$(2.5) \qquad \langle f, g \rangle = \int \overline{(\hat{f}(x))} \hat{g}(x) \, dx,$$

where $\overline{z}$ is the complex conjugate of $z$. Any continuous linear form $f \to u(f)$ on $H^s$ is of the form $u(f) = \langle g, f \rangle$ for some $g \in H^{-s}$, with $\|g\|_{H^{-s}} = \|u\|_{(H^s)^*}$, so that henceforth we shall identify the dual $(H^s)^*$ of $H^s$ with $H^{-s}$.

Fix some $s > 1/2$. We then have $H^s \subset C^0 \cap L^\infty$, so that $H^{-s}$ contains all bounded Radon measures on $\mathbb{R}$. In $H^s$, consider the set $H_-^s$ of functions with support in $]-\infty, 0]$, so that $f \in H_-^s$ if and only if $f(t) = 0$ for all $t > 0$. It is a closed subspace of $H^s$, so that the quotient space $H^s/H_-^s$ is a Hilbert space as well. This is the space we want:

$$(2.6) \qquad H = H^s/H_-^s.$$

To sum up, a real-valued function $f$ on $[0, \infty[$ belongs to $H$ if and only if it is the restriction to $[0, \infty[$ of some function in $H^s$, that is, if there is some function $\tilde{f} \in H^s$ (and hence defined on the whole real line) such that $\tilde{f}(t) = f(t)$ for all $t \ge 0$. The norm on $H$ is given by

$$\|f\|_H = \inf\{\|\tilde{f}\|_{H^s} | \tilde{f} \in H^s, \ \tilde{f}(t) = f(t) \ \forall t \ge 0\}$$

and the dual space $H^*$ by

$$H^* = \{g \in H^{-s} | \langle \tilde{f}, g \rangle = 0 \ \forall \tilde{f} \in H_-^s\}.$$

It follows that $H^*$ is the set of all distributions in $H^{-s}$ with support in $[0, \infty[$ and in particular, it contains all bounded Radon measures with support in $[0, \infty[$. $H$ inherits the property of being a Banach algebra from $H^s$.

2.2. *Bond dynamics.* From now on, it will be assumed that $p_t$ and $\bar{p}_t$ take values in $H$, so that the processes $p$ and $\bar{p}$ are $\mathcal{A}$-adapted and $H$-valued. Moreover, it will be assumed that $p$ and $\bar{p}$ are $\mathcal{A}$-progressively measurable. As in the finite-dimensional case, if $p$ (resp. $\bar{p}$) is $\mathcal{A}$-adapted and measurable, then it has an $\mathcal{A}$-progressively measurable modification. The reader who wants to avoid progressive measurability can therefore think of $p$ (resp. $\bar{p}$) as an $\mathcal{A}$-adapted measurable process.

We shall denote by $\mathcal{L}: [0, \infty[ \times H \to H$ the semigroup of left translations in $H$:

$$(2.7) \qquad (\mathcal{L}_a f)(s) = f(a+s),$$

where $a \geq 0$, $s \geq 0$ and $f \in H$. This is well defined since both $H^s$ and $H^s_-$ in (2.6) are invariant under left translations. One readily verifies that $\mathcal{L}$ is a strongly continuous contraction semigroup in $H$. Therefore (cf. Section 3, Chapter IX of [22]), it has an infinitesimal generator which we shall denote by $\partial$, with dense and invariant domain [the domain consists of all $f \in H$ such that $\lim_{\varepsilon \downarrow 0} \varepsilon^{-1}(\mathcal{L}_\varepsilon f - f)$ exists in $H$ and for such $f$ the limit is equal to $\partial f$], denoted by $\mathcal{D}(\partial)$. $\mathcal{D}(\partial)$ is a Hilbert space with norm

$$(2.8) \qquad \|f\|_{\mathcal{D}(\partial)} = (\|f\|_H^2 + \|\partial f\|_H^2)^{1/2}.$$

Volatilities are assumed to take values in the Hilbert space $\tilde{H}_0$ of all real-valued functions $F$ on $[0, \infty[$ such that $F = a + f$, for some $a \in \mathbb{R}$ and $f \in H$. The norm is given by

$$(2.9) \qquad \|F\|_{\tilde{H}_0} = (a^2 + \|f\|_H^2)^{1/2},$$

which is well defined since the decomposition of $F = a + f$, $a \in \mathbb{R}$ and $f \in H$, is unique. $\tilde{H}_0$ is a subset of continuous multiplication operators on $H$. In fact, since $H$ is a Banach algebra it follows that $\|Fh\|_H = C\|F\|_{\tilde{H}_0}\|h\|_H$, where $C > 0$ is independent of $F \in \tilde{H}_0$ and $h \in H$. We also introduce a Hilbert space $\tilde{H}_1$ of continuous multiplication operators on $\mathcal{D}(\partial)$. $\tilde{H}_1$ is the subspace of elements $F \in \tilde{H}_0$ with finite norm

$$(2.10) \qquad \|F\|_{\tilde{H}_1} = (a^2 + \|f\|_{\mathcal{D}(\partial)}^2)^{1/2},$$

where $F = a + f$, $a \in \mathbb{R}$ and $f \in \mathcal{D}(\partial)$. Finally let us define the left translation in $\tilde{H}_0$ by $(\mathcal{L}_a F)(s) = F(a+s)$, where $F \in \tilde{H}_0$, $a \geq 0$, $s \geq 0$. $\tilde{H}_1$ is the domain of the generator of $\mathcal{L}$, which we also denote $\partial$.

We shall assume that the bond dynamics are given by an equation of the following type. (Let $f_1, \ldots, f_n \in \tilde{H}_0$ (resp. $\tilde{H}_1$). Then $f_1 \cdots f_n \in \tilde{H}_0$ (resp. $\tilde{H}_1$) and when there is no risk for confusion, we shall also write $\mathcal{L}_a f_1 \cdots f_n$ instead of $\mathcal{L}_a(f_1 \cdots f_n)$. If moreover one $f_i \in H$ [resp. $\mathcal{D}(\partial)$], then $f_1 \cdots f_n \in H$ [resp. $\mathcal{D}(\partial)$].)

$$(2.11) \qquad \bar{p}_t = \mathcal{L}_t \bar{p}_0 + \int_0^t \mathcal{L}_{t-s}(m_s \bar{p}_s)\,ds + \int_0^t \sum_{i \in \mathbb{I}} \mathcal{L}_{t-s}(\sigma_s^i \bar{p}_s)\,dW_s^i,$$

for $t \in \mathbb{T}$, where $\sigma_t^i$, $i \in \mathbb{I}$, and $m_t$ are $\mathcal{A}$-progressively measurable $\tilde{H}_0$-valued processes and the $W^i$, $i \in \mathbb{I}$, are the already introduced standard Brownian motions. One must also take into account the boundary condition $B(t, t) = 1$, which in this context becomes

$$(2.12) \qquad \bar{p}_t(0) = \exp\left(-\int_0^t r(s)\, ds\right).$$

This can only be satisfied in general if

$$(2.13) \qquad \sigma_t^i(0) = 0 \qquad \text{for } i \in \mathbb{I}$$

and

$$(2.14) \qquad m_t(0) = 0.$$

When $\mathbb{I}$ is finite, then (2.11) gives the usual HJM equation (equation (9) of [7]) for $B$.

In this paper, the process $\bar{p}$ is given. So formula (2.11), which then defines $\sigma$ and $m$, can be considered as the decomposition of the real-valued semimartingale $t \mapsto \bar{p}_t(T - t) = \tilde{B}(t, T)$, describing the value of the zero-coupon bond with maturity $T$, for each fixed value of $T$. Alternatively, one may want to take $\sigma_t^i$ and $m_t$ as the parameters in the model, and derive $\bar{p}$ as the solution of a stochastic differential equation in $H$. Proceeding formally, (2.11) gives after differentiation

$$(2.15) \qquad \bar{p}_t = \bar{p}_0 + \int_0^t (\partial \bar{p}_s + \bar{p}_s m_s)\, ds + \int_0^t \bar{p}_s \sum_{i \in \mathbb{I}} \sigma_s^i\, dW_s^i.$$

A *mild solution* (cf. [4], Chapter 6, Section 1 for the case of deterministic $\sigma$) of (2.15) is an $\mathcal{A}$-progressively measurable $H$-valued process $\bar{p}$ satisfying the condition

$$(2.16) \qquad \int_0^{\bar{T}} \left( \|\bar{p}_t\|_H + \|\bar{p}_t m_t\|_H + \sum_{i \in \mathbb{I}} \|\bar{p}_t \sigma_t^i\|_H^2 \right) dt < \infty \qquad \text{a.s.}$$

and which satisfies (2.11). An $\mathcal{A}$-progressively measurable $H$-valued process $\bar{p}$ is a *strong solution* of (2.15) if condition (2.16) is satisfied and if $\bar{p}_t \in \mathcal{D}(\partial)$ a.s. for each $t \in \mathbb{T}$ and

$$(2.17) \qquad \int_0^{\bar{T}} \|\partial \bar{p}_t\|_H\, dt < \infty \qquad \text{a.s.}$$

We note that a strong solution of (2.15) is a semimartingale and it satisfies the evolution equation (2.11). However, the last term on the right-hand side of (2.11) is not in general the local martingale part. The aim of the following theorem is to ensure consistency in our model between the properties of $\bar{p}$ and those of $\sigma$ and $m$.

THEOREM 2.1. *If $\sigma^i$, $i \in \mathbb{I}$, and $m$ are given $\mathcal{A}$-progressively measurable $\tilde{H}_1$-valued processes, such that (2.13) and (2.14) are satisfied and such that*

$$\text{(2.18)} \qquad \int_0^{\bar{T}} \sum_{i \in \mathbb{I}} \|\sigma_t^i\|_{\tilde{H}_1}^2 \, dt < \infty \qquad a.s.$$

*and*

$$\text{(2.19)} \qquad \int_0^{\bar{T}} \|m_t\|_{\tilde{H}_1} \, dt < \infty \qquad a.s.$$

*and if $\bar{p}_0 \in H$ is given and satisfies [we use obvious functional notation such as $f > 0$ for $f \in H$, meaning $\forall s > 0 \; f(s) > 0$]*

$$\text{(2.20)} \qquad \bar{p}_0 \in \mathcal{D}(\partial), \qquad \bar{p}_0(0) = 1, \qquad \bar{p}_0 > 0,$$

*then (2.11) has, in the set of mild solutions of (2.15), a unique solution $\bar{p}$. This solution has the following properties: $\bar{p}$ is a strong solution of (2.15), $\bar{p}$ is strictly positive (i.e., $\forall t \in \mathbb{T}$, $\bar{p}_t > 0$), $t \mapsto \partial \bar{p}_t \in H$ is continuous a.s., the boundary condition*

$$\text{(2.21)} \qquad \bar{p}_t(0) = \exp\left( \int_0^t \frac{(\partial \bar{p}_s)(0)}{\bar{p}_s(0)} \, ds \right)$$

*is satisfied for each $t \in \mathbb{T}$ and an explicit expression of the solution is given by*

$$\text{(2.22)} \quad \bar{p}_t = \exp\left( \int_0^t \mathcal{L}_{t-s}\left( \left( m_s - \tfrac{1}{2} \sum_{i \in \mathbb{I}} (\sigma_s^i)^2 \right) ds + \sum_{i \in \mathbb{I}} \sigma_s^i \, dW_s^i \right) \right) \mathcal{L}_t \bar{p}_0.$$

*In particular $\bar{p}_t \in C^1([0, \infty[)$ a.s.*

So, given appropriate $\sigma^i$, $i \in \mathbb{I}$, and $m$, the mixed initial value and boundary value problem (2.11), (2.12) has a unique solution for any initial curve of zero-coupon bond prices satisfying (2.20). The proof of Theorem 2.1 is given in Section 5.

Under additional conditions on $\sigma^i$, $i \in \mathbb{I}$, and $m$, we are able to prove $L^p$-estimates of $\bar{p}$.

THEOREM 2.2. *If $\sigma^i$, $i \in \mathbb{I}$, and $m$ in Theorem 2.1 satisfy the following supplementary conditions: for each $a \in [1, \infty[$,*

$$\text{(2.23)} \quad E\left( \left( \int_0^{\bar{T}} \sum_{i \in \mathbb{I}} \|\sigma_t^i\|_{\tilde{H}_1}^2 \, dt \right)^a + \exp\left( a \int_0^{\bar{T}} \sum_{i \in \mathbb{I}} \|\sigma_t^i\|_{\tilde{H}_0}^2 \, dt \right) \right) < \infty$$

*and*

$$\text{(2.24)} \qquad E\left( \left( \int_0^{\bar{T}} \|m_t\|_{\tilde{H}_1} \, dt \right)^a + \exp\left( a \int_0^{\bar{T}} \|m_t\|_{\tilde{H}_0} \, dt \right) \right) < \infty,$$

*then the solution $\bar{p}$ in Theorem 2.1 has the following property: If $u \in [1, \infty[$, $q(t) = p_t / \mathcal{L}_t p_0$ and $\bar{q}(t) = \bar{p}_t / \mathcal{L}_t \bar{p}_0$, then $p, \bar{p} \in L^u(\Omega, P, L^\infty(\mathbb{T}, \mathcal{D}(\partial)))$ and $q, \bar{q}, 1/q, 1/\bar{q} \in L^u(\Omega, P, L^\infty(\mathbb{T}, \tilde{H}_1))$.*

We remind that, under the hypotheses of Theorem 2.1, $\bar{p}_t(0)$ satisfies (2.21), so it is the discount factor (2.12). Theorem 2.1 has the

COROLLARY 2.3. *Under the hypotheses of Theorem 2.2, if $\alpha \in \mathbb{R}$, then the discount factor satisfies*

$$E\left(\sup_{t \in \mathbb{T}} (\bar{p}_t(0))^\alpha\right) < \infty.$$

2.3. *Portfolios.* The linear functionals in $H^*$ will be interpreted as bond portfolios. More precisely, a portfolio is an $H^*$-valued $\mathcal{A}$-progressively measurable process $\theta$ defined on $\mathbb{T}$. Its value at time $t$ is

(2.25) $$V(t, \theta) = \langle \theta_t, p_t \rangle$$

and its discounted value is

(2.26) $$\bar{V}(t, \theta) = \langle \theta_t, \bar{p}_t \rangle.$$

EXAMPLE 2.4. (i) A portfolio consisting of one single zero-coupon bond with a fixed *time of maturity* $T$, $T \geq \bar{T}$, is represented by $\theta$, where $\theta_t = \delta_{T-t} \in H^*$, the Dirac mass with support at $T - t$, where $t \in \mathbb{R}$. Note that when $t$ increases its support moves to the left toward the origin, which also can be expressed by $\theta_t(s) = \theta_0(s + t)$, for $s \geq 0$. Its value at time $t$ is $p_t(T - t)$.

(ii) A portfolio $\theta$ consisting of one single zero-coupon bond with a fixed *time of maturity* $T$, $0 \leq T < \bar{T}$. Then $\theta_t = \delta_{T-t} \in H^*$, for $t \leq T$ and $\theta_t = 0$, for $T < t \leq \bar{T}$. Its value at time $t \leq T$ is $p_t(T - t)$ and its value at time $t > T$ is zero.

(iii) $\theta$ given by $\theta_t = \delta_S \in H^*$, the Dirac mass with fixed support at $S$, represents a portfolio which consists at any time of a single zero-coupon bond with *time to maturity* $S$; note that it has to be constantly readjusted to keep the time to maturity constant, and that its value at time $t$ is $p_t(S)$.

As usual, a portfolio will be called *self-financing* if at any time, the change in its value is due to changes in market prices, and not to any redistribution of the portfolio, that is,

(2.27) $$\bar{V}(t, \theta) = \bar{V}(0, \theta) + \bar{G}(t, \theta),$$

where $\bar{G}(t, \theta)$ represents the discounted gains in the time interval $[0, t[$. We shall find the expression of $\bar{G}(t, \theta)$. We remind that the subspace of elements $f$ of $H^*$ with support not containing $0$ is dense in $H^*$. Suppose that the portfolio is already defined up to time $t$ and that $\theta_t$ contains no zero-coupon bonds of time to maturity smaller than some $A > 0$, that is, $\theta_t$ has no support in $[0, A[$. At $t$ let the portfolio evolve itself without any trading until $t + \varepsilon$, where $0 < \varepsilon < A$. Then $\theta_{t+\varepsilon}$ is given by $\theta_{t+\varepsilon}(s) = \theta_t(s + \varepsilon)$, for $s \geq 0$. At $t + \varepsilon$, the discounted value of the portfolio is $\bar{V}(t + \varepsilon, \theta) = \langle \theta_{t+\varepsilon}, \bar{p}_{t+\varepsilon} \rangle = \int_A^\infty \theta_t(s) \bar{p}_{t+\varepsilon}(s - \varepsilon) \, ds$. We can now differentiate in $\varepsilon$. Using (2.11) and (2.27) and taking the limits $\varepsilon \to 0$ and then $A \to 0$ we obtain:

$$(2.28) \qquad d\bar{G}(t, \theta) = \langle \theta_t, \bar{p}_t m_t \rangle \, dt + \sum_{i \in \mathbb{I}} \langle \theta_t, \bar{p}_t \sigma_t^i \rangle \, dW_t^i.$$

We now take $\bar{G}(0, \theta) = 0$ and this expression, in case it makes sense, as the definition of the discounted gains for an arbitrary portfolio $\theta$.

To formalize this idea, we need to define appropriately the space of admissible portfolios. Given the process $\bar{p}$, an *admissible portfolio* is an $H^*$-valued $\mathcal{A}$-progressively measurable process $\theta$ such that

$$
\begin{aligned}
(2.29) \qquad \|\theta\|_{\mathsf{P}}^2 = E \bigg( \int_0^{\bar{T}} (\|\theta_t\|_{H^*}^2 &+ \|\sigma_t^* \theta_t \bar{p}_t\|_{H^*}^2) \, dt \\
&+ \bigg( \int_0^{\bar{T}} |\langle \theta_t, \bar{p}_t m_t \rangle| \, dt \bigg)^2 \bigg) < \infty,
\end{aligned}
$$

where we have used the *notation*

$$(2.30) \qquad \|\sigma_t^* \theta_t \bar{p}_t\|_{H^*}^2 = \sum_{i \in \mathbb{I}} (\langle \theta_t, \bar{p}_t \sigma_t^i \rangle)^2.$$

For the mathematically minded reader this notation will be given a meaning in Section 5. The set of all admissible portfolios is Banach space $\mathsf{P}$ and the subset of all admissible self-financing portfolios is denoted by $\mathsf{P}_{\mathrm{sf}}$. The discounted gains process for a portfolio in $\mathsf{P}$ is a continuous square-integrable process:

PROPOSITION 2.5. *Assume that $\bar{p}_0$, $m$ and $\sigma$ are as in Theorem* 2.1. *If $\theta \in \mathsf{P}$, then $\bar{G}(\cdot, \theta)$ is continuous a.s. and $E(\sup_{t \in \mathbb{T}} (\bar{G}(t, \theta))^2) < \infty$.*

EXAMPLE 2.6. (i) The portfolio of Example 2.4(i) is self-financing and the portfolios of Example 2.4(ii) and (iii) are not self-financing.

(ii) We define a self-financing portfolio $\theta$ of zero-coupon bonds with constant *time to maturity* $S$. Let $\theta$ be given by $\theta_t = x(t) \delta_S$, where

$$(2.31) \qquad x(t) = x(0) \exp \bigg( \int_0^t f_s(S) \, ds \bigg)$$

and $f_t(S)$ is given by (2.1). That $\theta$ is self-financing is readily established by observing that in this case $x(t) \bar{p}_t(S) = \bar{V}(t, \theta)$, $\bar{V}(t, \theta) = \bar{V}(0, \theta) + \int_0^t \bar{V}(s, \theta) \times$ $(m_s(S)\,ds + \sum_{i\in\mathbb{I}}\sigma_s^i(S)\,dW_s^i)$ and by applying Itô's lemma to $x(t) = \bar{V}(t,\theta)/\bar{p}_t(S)$; cf. [19].

We note that $x(t) = V(t,\theta)/p_t(S)$ is the wealth at time $t$ expressed in units of zero-coupon bonds of time to maturity $S$. According to (2.31), the self-financing portfolio $\theta$ is then given by the initial number $x(0)$ of bonds and by the growth rate $f(S)$ of $x$. So this is as a money account, except that here we count in zero-coupon bonds of time to maturity $S$.

In particular, if $S = 0$, then the equality $x(t) = V(t,\theta)$, the definition (2.1) of $r$ and the definition (2.31) show that $\theta$ can be assimilated to money at a usual bank account with spot rate $r$, see [2].

REMARK 2.7 (Roll-Overs). (i) Let $S \geq 0$, $x(0) = 1$ and the portfolio $\theta$ be as in (ii) of Example 2.6. Of course $X = V(\bar{T},\theta)$ is then an attainable interest rate derivative, for which $\theta$ is a replicating portfolio. We name this derivative a *Roll-Over* or more precisely an $S$-Roll-Over to specify the time to maturity of the underlying zero-coupon bond. Let $\tilde{p}_t(S)$ be the discounted price of an $S$-Roll-Over at time $t$. Then $\tilde{p}_0(S) = p_0(S)$ by definition and the price dynamics of Roll-Overs is simply given by

$$(2.32) \qquad \tilde{p}_t = p_0 + \int_0^t \tilde{p}_s m_s\,ds + \int_0^t \tilde{p}_s \sum_{i\in\mathbb{I}} \sigma_s^i\,dW_s^i,$$

$t \in \mathbb{T}$, which solution $\tilde{p}$ is given by

$$(2.33) \qquad \tilde{p}_t(S) = \bar{p}_t(S)\exp\left(\int_0^t f_s(S)\,ds\right), \qquad S \geq 0.$$

An $S$-Roll-Over can be denounced at time $t$, with a notice of $S$ time units and it will then pay $x(t)$ units of account at time $t + S$.

(ii) Zero-coupon bonds do not in general permit self-financing buy-and-hold portfolios, that is, constant portfolios. However, Roll-Overs do, since a constant portfolio of Roll-Overs is always self-financing. Mathematically, this can be thought of as changing from a fixed frame to a moving frame for expressing a self-financed discounted wealth process in terms of coordinates, that is the portfolio. To be more precise let us consider a technically simple case. Let $\sigma$ be nondegenerated in the sense that the linear span of the set $\{\sigma_t^i | i \in \mathbb{I}\}$ is dense a.s. in $\tilde{H}_0$ for every $t \in \mathbb{T}$. Let the initial price satisfy $\sup_{t\in\mathbb{T}}\sup_{s\geq 0} p_0(s)/p_0(t+s) < \infty$ and let the hypotheses of Theorem 2.2 be satisfied. Then a self-financing portfolio $\theta \in \mathsf{P}_{sf}$ is the unique replicating portfolio in $\theta \in \mathsf{P}_{sf}$ of $V(\bar{T},\theta)$. Moreover, there is a unique $\eta \in \mathsf{P}$ such that $\langle \theta_t, \bar{p}_t \rangle = \langle \eta_t, \tilde{p}_t \rangle$, for all $t \in \mathbb{T}$. The coordinates of the self-financed discounted wealth process $V(\cdot,\theta)$ with respect to the moving frame is $\eta$. In particular, an $S$-Roll-Over is given by the constant portfolio $\eta$, where $\eta_t = \delta_S$.

We next set up an arbitrage-free market by postulating a market-price of risk relation between $m$ and $\sigma$.

CONDITION A. *There exists a family $\{\Gamma^i | i \in \mathbb{I}\}$ of real-valued $\mathcal{A}$-progressively measurable processes such that*

$$(2.34) \qquad m_t + \sum_{i \in \mathbb{I}} \Gamma_t^i \sigma_t^i = 0$$

*and*

$$(2.35) \qquad E\left( \exp\left( a \int_0^{\bar{T}} \sum_{i \in \mathbb{I}} |\Gamma_t^i|^2 \, dt \right) \right) < \infty \qquad \forall \, a \geq 0.$$

Condition (2.34) is similar to a standard no-arbitrage condition in finite dimension and we refer to Remark 5.1 for further motivation in the infinite-dimensional case. Inequality (2.35) permits the use of Novikov's criteria (cf. [18], Chapter VIII, Proposition 1.15). When Condition A is satisfied, (2.28) for the discounted gains of a portfolio $\theta$ becomes

$$(2.36) \qquad d\bar{G}(t, \theta) = \sum_{i \in \mathbb{I}} \langle \theta_t, \bar{p}_t \sigma_t^i \rangle (-\Gamma_t^i \, dt + dW_t^i).$$

The following result shows how to obtain a martingale measure in the general case of Condition A. Introduce the notation

$$(2.37) \qquad \xi_t = \exp\left( -\tfrac{1}{2} \int_0^t \sum_{i \in \mathbb{I}} (\Gamma_s^i)^2 \, ds + \int_0^t \sum_{i \in \mathbb{I}} \Gamma_s^i \, dW_s^i \right),$$

where $t \in \mathbb{T}$.

THEOREM 2.8. *If (2.35) is satisfied, then $\xi$ is a martingale with respect to $(P, \mathcal{A})$ and $\sup_{t \in \mathbb{T}} \xi_t^\alpha \in L^1(\Omega, P)$ for each $\alpha \in \mathbb{R}$. The measure $Q$, defined by*

$$dQ = \xi_{\bar{T}} \, dP,$$

*is equivalent to $P$ on $\mathcal{F}_{\bar{T}}$ and $t \mapsto \bar{W}_t^i = W_t^i - \int_0^t \Gamma_s^i \, ds$, $t \in \mathbb{T}$, $i \in \mathbb{I}$, are independent Wiener process with respect to $(Q, \mathcal{A})$. (The Girsanov formula holds.)*

The expected value of a random variable $X$ with respect to $Q$ is denoted $E_Q(X)$ and $E_Q(X) = E(\xi_{\bar{T}} X)$.

Proposition 2.5 and Theorem 2.8 have the

COROLLARY 2.9. *Assume that $\bar{p}_0$ and $\sigma$ are as in Theorem 2.1 and assume that Condition A is satisfied. Then all conditions of Theorem 2.1 are satisfied and if $\theta \in \mathsf{P}$, then $\bar{G}(\cdot, \theta)$ is continuous a.s., $E(\sup_{t \in \mathbb{T}} (\bar{G}(t, \theta))^2) < \infty$ and $\bar{G}(\cdot, \theta)$ is a $(Q, \mathcal{A})$-martingale.*

By an arbitrage-free market, we mean as usually, that there does not exist a self-financing dynamical portfolio $\theta \in \mathsf{P}_{\mathrm{sf}}$ such that $V(0, \theta) = 0$, $V(\bar{T}, \theta) \geq 0$ and $P(V(\bar{T}, \theta) > 0) > 0$. The following result shows that the market is arbitrage-free:

COROLLARY 2.10.    *Assume that $\bar{p}_0$ and $\sigma$ are as in Theorem 2.1 and assume that Condition A is satisfied. If $\theta \in \mathsf{P}_{\mathrm{sf}}$, its discounted price $\bar{V}(\cdot, \theta)$ is a $(Q, \mathcal{A})$-martingale and $E(\sup_{t \in \mathbb{T}} (\bar{V}(t, \theta))^2) < \infty$. In particular the market is arbitrage-free.*

**3. The optimal portfolio problem.**  The investor is characterized by his utility $u(\bar{w}_{\bar{T}})$, where $\bar{w}_{\bar{T}}$ is terminal wealth, discounted to $t = 0$. Given the initial wealth $x$, denote by $\mathcal{C}(x)$ the set of all admissible self-financing portfolios starting from $x$:

$$\mathcal{C}(x) = \{\theta \in \mathsf{P}_{\mathrm{sf}} | \bar{V}(0, \theta) = x\}.$$

The investor's optimization problem is, for a given initial wealth $K_0$, to find a solution $\hat{\theta} \in \mathcal{C}(K_0)$ of

$$(3.1) \qquad E(u(\bar{V}(\bar{T}, \hat{\theta}))) = \sup_{\theta \in \mathcal{C}(K_0)} E(u(\bar{V}(\bar{T}, \theta))).$$

In the following, the utility function is allowed to take the value $-\infty$, so $u : \mathbb{R} \to \mathbb{R} \cup \{-\infty\}$. Throughout this section, we make the following Inada-type assumptions:

CONDITION B.

(a) $u : \mathbb{R} \to \mathbb{R} \cup \{-\infty\}$ *is strictly concave, upper semi-continuous and finite on an interval* $]\underline{x}, \infty[$, *with* $\underline{x} \leq 0$ *(the value* $\underline{x} = -\infty$ *is allowed ).*

(b) *$u$ is $C^1$ on* $]\underline{x}, \infty[$ *and* $u'(x) \to \infty$ *when* $x \to \underline{x}$ *in* $]\underline{x}, \infty[$.

(c) *There exists some $q > 0$ such that*

$$(3.2) \qquad \liminf_{x \downarrow \underline{x}} (1 + |x|)^{-q} u'(x) > 0$$

*and such that, if $u' > 0$ on* $]\underline{x}, \infty[$, *then*

$$(3.3) \qquad \limsup_{x \to \infty} x^q u'(x) < \infty$$

*and if $u'$ takes the value zero, then*

$$(3.4) \qquad \limsup_{x \to \infty} x^{-q} u'(x) < 0.$$

REMARK 3.1. (i) If $u$ satisfies Condition B, then $v$ obtained by an affine transformation, $v(x) = \alpha u(ax + b) - \beta$, $\alpha, a > 0$, $\beta \in \mathbb{R}$, $b \geq \underline{x}$, also satisfies Condition B. Usual utility functions, such as exponential $u(x) = -e^{-x}$, quadratic $u(x) = -x^2/2$, power $u(x) = x^a/a$, $x > 0$, $a < 1$ and $a \neq 0$ and logarithmic $u(x) = \ln x$, $x > 0$, satisfy Condition B. Others, like HARA, are obtained by affine transformations.

(ii) Strictly negative wealth is admitted when $\underline{x} < 0$ and additional constraints such as positivity are not included in the present theory. However, positivity of wealth is obviously satisfied for all utility functions with $\underline{x} = 0$, such as $u(x) = x^a/a$, $a < 1$ and $a \neq 0$ and logarithmic $u(x) = \ln x$.

It follows that $u'$ restricted to $]\underline{x}, \infty[$ has a strictly decreasing continuous inverse $\varphi$, that is, a map such that $(\varphi \circ u')(x) = x$ for $x \in ]\underline{x}, \infty[$. The domain of $\varphi$ is $I = u'(]\underline{x}, \infty[)$. Condition B has an equivalent formulation in terms of $\varphi$:

LEMMA 3.2. *If $u$ satisfies Condition B, then:*

(i) *If $u' > 0$ on $]\underline{x}, \infty[$, then $I = ]0, \infty[$ and for some $C, p > 0$,*

$$(3.5) \qquad |\varphi(x)| \leq C(x^p + x^{-p}),$$

*for all $x > 0$.*

(ii) *If $u'$ takes the value zero in $]\underline{x}, \infty[$, then $I = \mathbb{R}$ and for some $C, p > 0$,*

$$(3.6) \qquad |\varphi(x)| \leq C(1 + |x|)^p,$$

*for all $x \in \mathbb{R}$.*

*Conversely, if $I = ]0, \infty[$ (resp. $I = \mathbb{R}$), $\underline{x} \in [-\infty, 0]$, and $\varphi : I \to ]\underline{x}, \infty[$ satisfying (3.5) [resp. (3.6)] is a strictly decreasing continuous surjection with inverse $g$, $a \in \mathbb{R}$, $a \in ]\underline{x}, \infty[$ and if $u(x) = a + \int_{x_0}^{x} g(y) \, dy$, for $x \in ]\underline{x}, \infty[$, $u(\underline{x}) = \lim_{x \downarrow \underline{x}} u(x)$, $u(x) = -\infty$, for $x < \underline{x}$, then $u$ satisfies Condition B.*

We shall next give existence results of optimal portfolios. In order to construct solutions of the optimization problem (3.1), we first solve a related problem of optimal terminal discounted wealth at time $\bar{T}$, which gives candidates of optimal terminal discounted wealths, and second, we construct, for certain of these candidates, a hedging portfolio, which then is a solution of the optimal portfolio problem (3.1). The construction of terminal discounted wealths is general and only requires that Conditions A and B are satisfied. For the construction of hedging portfolios, we separate the case of a finite number of random sources, that is, $\mathbb{I} = \{1, \ldots, \bar{m}\}$ (Theorem 3.6) and the case of infinitely many random sources $\mathbb{I} = \mathbb{N}^*$ (Theorem 3.8). In the case of $\mathbb{I}$ finite, general stochastic volatilities being nondegenerated according to a certain condition are considered. In the case of $\mathbb{I}$ infinite, we only give results for deterministic $\sigma$, but which can be degenerated.

If $X$ is the terminal discounted wealth for a self-financing strategy in $\mathcal{C}(K_0)$, then due to Corollary 2.9 $K_0 = E(\xi_{\bar{T}}X)$. We shall employ dual techniques to find candidates of the optimal $X$; cf. [17].

THEOREM 3.3. *Let $u$ satisfy Condition* B *and let $\Gamma$ satisfy condition* (2.35). *If $K_0 \in ]\underline{x}, \infty[$, then there exists a unique $\hat{X} \in L^2(\Omega, P, \mathcal{F}_{\bar{T}})$ such that $K_0 = E(\xi_{\bar{T}}\hat{X})$ and*

$$(3.7) \qquad E(u(\hat{X})) = \sup_{K_0 = E(\xi_{\bar{T}}X)} E(u(X)).$$

*Moreover, $\hat{X} \in L^p(\Omega, P)$ for each $p \in [1, \infty[$ and there is a unique $\lambda \in I$ such that $\hat{X} = \varphi(\lambda \xi_{\bar{T}})$.*

Now, if $\hat{\theta}$ hedges $\hat{X}$, then $\hat{\theta}$ is an optimal portfolio. More precisely, we have:

COROLLARY 3.4. *Let $m$ and $\sigma$ satisfy the hypotheses of Theorem 2.1 and also be such that there exists a $\Gamma$ with the following properties: $\Gamma$ satisfies Condition* A *and $\hat{X}$, given by Theorem 3.3, satisfies $\hat{X} = \bar{V}(\bar{T}, \hat{\theta})$ for some $\hat{\theta} \in \mathsf{P}_{\mathrm{sf}}$. Then $\hat{\theta}$ is a solution of the optimal portfolio problem* (3.1).

REMARK 3.5. Instead of optimizing in (3.1) the expected utility of the discounted terminal wealth $\bar{V}(\bar{T}, \theta)$, one can choose to optimize that of the terminal wealth $V(\bar{T}, \theta)$. This leads to a similar result as that of Theorem 3.3. Using that $\bar{p}_t(0)$ is the discount factor, we obtain that the optimal terminal wealth is $\hat{Z} = \varphi(\bar{p}_{\bar{T}}(0)\lambda\xi_{\bar{T}})$ and that $\lambda$ is given by $E(\xi_{\bar{T}}\bar{p}_{\bar{T}}(0)\hat{Z}) = K_0$.

3.1. *The case $\mathbb{I} = \{1, \ldots, \bar{m}\}$.* Here we assume that a.s. the set of volatilities $\{\sigma_t^1, \ldots, \sigma_t^{\bar{m}}\}$ is linearly independent in $\tilde{H}_0$ for each $t \in \mathbb{T}$. Since $p_0 > 0$ and $p_0 \in H$, this is equivalent to the a.s. linear independence of $\{\sigma_t^1 \mathcal{L}_t p_0, \ldots, \sigma_t^{\bar{m}} \mathcal{L}_t p_0\}$ in $H$, for each $t$. Consider the $\bar{m} \times \bar{m}$ matrix $A(t)$ with elements

$$A(t)_{ij} = (\sigma_t^i \mathcal{L}_t p_0, \sigma_t^j \mathcal{L}_t p_0)_H$$

(beware that we are using the scalar product in $H$ and not in $L^2$):

THEOREM 3.6. *Let $p_0 \in \mathcal{D}(\partial)$, $p_0(0) = 1$ and $p_0 > 0$, let $\sigma \neq 0$ satisfy conditions* (2.13) *and* (2.23) *and let Conditions* A *and* B *be satisfied. Assume that there exists an adapted process $k > 0$, such that for each $q \geq 1$ we have $E(\sup_{t \in \mathbb{T}} k_t^q) < \infty$ and, for each $x \in \mathbb{R}^{\bar{m}}$ and $t \in \mathbb{T}$:*

$$(3.8) \qquad (x, A(t)x)_{\mathbb{R}^{\bar{m}}} k_t \geq \left( \sum_{i \in \mathbb{I}} \|\sigma_t^i \mathcal{L}_t p_0\|_H^2 \right)^{1/2} \|x\|_{\mathbb{R}^{\bar{m}}}^2 \qquad a.s.$$

*If $K_0 \in \, ]\underline{x}, \infty[$, then problem (3.1) has a solution $\hat{\theta}$.*

We note that condition (3.8) only involves prices at time 0 and the volatilities. We also note that the optimal portfolio is never unique since one can always add a nontrivial portfolio $\theta'$ such that the linear span of the set $\{\sigma_t^1 \bar{p}_t, \ldots, \sigma_t^{\bar{m}} \bar{p}_t\}$ is in the kernel of $\theta'_t$.

REMARK 3.7. Due to the nonuniqueness of the optimal portfolio $\hat{\theta}$ in Theorem 3.6, it can be realized using different numbers of bonds:

1. One can always choose an optimal portfolio $\hat{\theta}$ such that $\hat{\theta}_t$ consists of at most $1 + \bar{m}$ zero-coupon bonds at every time $t$. This can be seen by a heuristic argument. Since, for every $t \geq 0$, the set of continuous functions $\{\sigma_t^1 \bar{p}_t, \ldots, \sigma_t^{\bar{m}} \bar{p}_t\}$ is linearly independent a.s., there exists positive $\mathcal{F}_t$-measurable finite random variables $S_t^j$ such that $0 < S_t^1 < \cdots < S_t^{\bar{m}}$ and such that the vectors $v_t^j = (\sigma_t^1(S_t^j)\bar{p}_t(S_t^j), \ldots, \sigma_t^{\bar{m}}(S_t^j)\bar{p}_t(S_t^j))$, $1 \leq j \leq \bar{m}$, are linearly independent a.s. Let $\theta_t = \sum_{1 \leq j \leq \bar{m}} a_t^j \delta_{S_t^j}$, where $a_t^j$ are real $\mathcal{F}_t$-measurable random variables. The equations $\langle \theta_t, \bar{p}_t \sigma_t^i \rangle = y_i(t)$, $1 \leq i \leq \bar{m}$, where $y_i(t)$ is given by (5.31), then have a unique solution $a_t$. So at time $t$ it is enough to use bonds with time to maturity $0 = S_t^0 < S_t^1 < \cdots < S_t^{\bar{m}}$ to realize an optimal portfolio $\hat{\theta}$. The number of bonds with time to maturity $0 = S_t^0$ is adjusted to obtain a self-financing portfolio.

2. Alternatively to zero-coupon bonds, one can also use $\bar{m} + 1$ coupon bonds or Roll-Overs to realize an optimal portfolio.

3. For certain volatility structures, one can even use any $\bar{m}$ given Roll-Overs or $\bar{m}$ given marketed coupon bonds (supposed to have distinct times of maturity, each exceeding $\bar{T}$) to realize an optimal portfolio. In particular, this is the case if the above vectors $v_t^j$, $1 \leq j \leq \bar{m}$, are linearly independent for every sequence $0 < S_t^1 < \cdots < S_t^{\bar{m}}$.

3.2. *The case of deterministic $\sigma$ and $\Gamma$.* Condition (3.8) cannot hold in the infinite case, $\mathbb{I} = \mathbb{N}^*$. In fact this is a consequence of that $\int_0^{\bar{T}} \sum_{i \in \mathbb{I}} \|\sigma_t^i\|_{\bar{H}_0}^2 \, dt < \infty$ a.s., as explained in Remark 5.3. When $\sigma$ and $\Gamma$ are deterministic, we can give another, weaker condition, which only involves $\sigma$, $\Gamma$ and the zero-coupon bond prices at time zero $\bar{p}_0$. This will give us a result which will hold for the infinite case as well. Properties of the inverse $\varphi$ of the derivative of the utility function $u$, satisfying Condition B, were given in Lemma 3.2. For simplicity we shall need one more property, which we impose directly as a condition on $\varphi$. We keep in mind that $\varphi' < 0$, since $u$ is strictly concave.

CONDITION C. *Let Condition B be satisfied, assume that $u$ is $C^2$ on $]\underline{x}, \infty[$ and assume that there exist $C, p > 0$ such that:*

(a) *If $u' > 0$ on $]\underline{x}, \infty[$, then*

$$(3.9) \qquad |x\varphi'(x)| \le C(x^p + x^{-p}),$$

*for all $x > 0$.*

(b) *If $u'$ takes the value zero in $]\underline{x}, \infty[$, then*

$$(3.10) \qquad |x\varphi'(x)| \le C(1 + |x|)^p,$$

*for all $x \in \mathbb{R}$.*

We note that Condition B implies Condition C if $u'$ is homogeneous. Condition C is satisfied for the utility functions in Remark 3.1.

THEOREM 3.8. *Let $\sigma$ and $m$ be deterministic, while $\mathbb{I}$ is finite or infinite. Let $p_0 \in \mathcal{D}(\partial)$, $p_0(0) = 1$ and $p_0 > 0$, let $\sigma$ satisfy conditions (2.13) and (2.18) and let Conditions A and C be satisfied. Assume that there exists a (deterministic) $H^*$-valued function $\gamma \in L^2(\mathbb{T}, H^*)$ such that*

$$(3.11) \qquad \langle \gamma_t, \sigma_t^i \mathcal{L}_t p_0 \rangle = \Gamma_t^i,$$

*for each $i \in \mathbb{I}$ and $t \in \mathbb{T}$. If $K_0 \in ]\underline{x}, \infty[$, then problem (3.1) has a solution $\hat{\theta}$.*

As explained in Remark 5.4, condition (3.11) can be satisfied in highly incomplete markets. In the situation of Theorem 3.8, we can derive an explicit expression of an optimal portfolio. We use the notation $\bar{q}(t) = \bar{p}_t / \mathcal{L}_t \bar{p}_0$ of Theorem 2.2.

COROLLARY 3.9. *Under the hypotheses of Theorem 3.8, an optimal portfolio is given by $\hat{\theta} = \theta^0 + \theta^1$, where $\theta^0, \theta^1 \in \mathsf{P}$, $\theta^0 = a_t \delta_0$, $\theta_t^1 = b_t(\bar{q}(t))^{-1} \gamma_t$. The coefficients $a$ and $b$ are real-valued $\mathcal{A}$-progressively measurable processes given by*

$$(3.12) \qquad b_t = E_Q(\lambda \xi_{\bar{T}} \varphi'(\lambda \xi_{\bar{T}}) | \mathcal{F}_t)$$

*and*

$$(3.13) \qquad a_t = (\bar{p}_t(0))^{-1}(Y(t) - b_t \langle \gamma_t, \mathcal{L}_t p_0 \rangle),$$

*$t \in \mathbb{T}$, where $Y(t) = E_Q(\varphi(\lambda \xi_{\bar{T}}) | \mathcal{F}_t)$ and $\lambda \in I$ is unique. The discounted price of the portfolio $\hat{\theta}$ is given by $\bar{V}(t, \hat{\theta}) = Y(t)$, $t \in \mathbb{T}$. Moreover, $\langle \theta^0, \sigma_t^i \bar{p}_t \rangle = 0$ and*

$$(3.14) \qquad \langle \theta^1, \sigma_t^i \bar{p}_t \rangle = E_Q(\lambda \xi_{\bar{T}} \varphi'(\lambda \xi_{\bar{T}}) | \mathcal{F}_t) \Gamma_t^i,$$

*$i \in \mathbb{I}$ and $t \in \mathbb{T}$.*

The proofs of Theorem 3.8 and Corollary 3.9 are based on a Clark–Ocone like representation of the optimal terminal discounted wealth (see Lemma A.5). Alternatively, the explicit expressions in Corollary 3.9 can be obtained by a Hamilton–Jacobi–Bellman approach (see Remark 5.5). This corollary has an important consequence, since it leads directly to mutual fund theorems. We shall state a version only involving self-financing portfolios.

THEOREM 3.10. *Under the hypotheses of Theorem* 3.8, *there exists a self-financing portfolio* $\Theta \in \mathsf{P}_{\mathrm{sf}}$, *with the following properties:*

(i) *The initial value of* $\Theta$ *is 1 euro, that is,* $\langle \Theta_0, \bar{p}_0 \rangle = 1$ *and the value at each time* $t \in \mathbb{T}$ *is strictly positive, that is,* $\langle \Theta_t, \bar{p}_t \rangle > 0$.

(ii) *For each given utility function* $u$ *satisfying Condition* C *and each initial wealth* $K_0 \in ]\underline{x}, \infty[$, *there exist two real-valued processes* $x$ *and* $y$ *such that if* $\hat{\theta}_t = x_t \delta_0 + y_t \Theta_t$, *then* $\hat{\theta}$ *is an optimal self financing portfolio for* $u$, *that is a solution of problem* (3.1).

**4. Examples of closed-form solutions.** In this section we shall give, in the situation of Corollary 3.9, examples of solutions of problem (3.1), for certain utility functions $u$. In particular, condition (3.11) is satisfied, so $\sigma$ and $\Gamma$ are deterministic.

According to Corollary 3.9, the final optimal wealth is $Y(\bar{T}) = \varphi(\lambda \xi_{\bar{T}})$ and the optimal discounted wealth process $Y$ is given by $Y(t) = E_Q(\varphi(\lambda \xi_{\bar{T}})|\mathcal{F}_t)$. The initial wealth $Y(0) = K_0$ determines $\lambda$. We introduce the optimal utility $U$ as a function of discounted wealth $w$ at time $t \in \mathbb{T}$,

$$(4.1) \qquad U(t, w) = E(u(Y(\bar{T}))|Y(t) = w).$$

We recall that $(\bar{p}_t)^{-1} \mathcal{L}_t \bar{p}_0 \in \tilde{H}_1$ a.s. and $(\bar{p}_t(0))^{-1} \in \mathbb{R}$ a.s. which is a particular case of Theorem 2.2 and Corollary 2.3.

EXAMPLE 4.1. *Quadratic utility.* The utility function is

$$(4.2) \qquad u(x) = \mu x - x^2/2,$$

where $\mu \in \mathbb{R}$ is given. Condition B is satisfied with $\underline{x} = -\infty$ and $0 < q < 1$. We have $I = \mathbb{R}$ and $\varphi(x) = -x + \mu$ and Condition C is satisfied with $p \geq 2$. The $P$-martingale $\xi$ can be written $\xi_t = \eta_t \exp(\int_0^t \sum_{i \in \mathbb{I}} (\Gamma_s^i)^2 \, ds)$, where $\eta_t = \exp(-\frac{1}{2} \int_0^t \sum_{i \in \mathbb{I}} (\Gamma_s^i)^2 \, ds + \int_0^t \sum_{i \in \mathbb{I}} \Gamma_s^i \, d\bar{W}_s^i)$ defines a $Q$-martingale $\eta$ (see Theorem 2.8). Since

$$K_0 = E_Q(\varphi(\lambda \xi_{\bar{T}})) = -\lambda E_Q(\xi_{\bar{T}}) + \mu = -\lambda \exp\left(\int_0^{\bar{T}} \sum_{i \in \mathbb{I}} (\Gamma_s^i)^2 ds\right) + \mu,$$

it follows that

$$(4.3) \qquad \lambda = (\mu - K_0) \exp\left(-\int_0^{\bar{T}} \sum_{i \in \mathbb{I}} (\Gamma_s^i)^2 \, ds\right).$$

The optimal discounted wealth process $Y$ is then given by $Y(t) = E_Q(\varphi(\lambda \xi_{\bar{T}})|\mathcal{F}_t) = -\lambda \exp(\int_0^{\bar{T}} \sum_{i \in \mathbb{I}} (\Gamma_s^i)^2 \, ds) E_Q(\eta_{\bar{T}}|\mathcal{F}_t) + \mu$, so

$$(4.4) \qquad Y(t) = \mu + (K_0 - \mu) \exp\left(\int_0^t \sum_{i \in \mathbb{I}} (\Gamma_s^i \, d\bar{W}_s^i - \tfrac{1}{2}(\Gamma_s^i)^2 \, ds)\right).$$

For given $\mu, w \in \mathbb{R}$ and $t \in \mathbb{T}$, formula (4.1) leads to the optimal utility

$$(4.5) \qquad \begin{aligned} U(t, w) &= (-\tfrac{1}{2}w^2 + \mu w) \exp\left(-\int_t^{\bar{T}} \sum_{i \in \mathbb{I}} (\Gamma_s^i)^2 \, ds\right) \\ &\quad + \tfrac{1}{2}\mu^2 \left(1 - \exp\left(-\int_t^{\bar{T}} \sum_{i \in \mathbb{I}} (\Gamma_s^i)^2 \, ds\right)\right). \end{aligned}$$

One finds first $b_t = Y(t) - \mu$ and then

$$(4.6) \qquad a_t = (\bar{p}_t(0))^{-1}(Y(t) - (Y(t) - \mu)\langle \gamma_t, \mathcal{L}_t p_0 \rangle).$$

An optimal portfolio is given by $\hat{\theta} = \theta^0 + \theta^1$, where $\theta_t^0 = a_t \delta_0$ and

$$(4.7) \qquad \theta_t^1 = (Y(t) - \mu)\gamma_t(\bar{p}_t)^{-1}\mathcal{L}_t p_0.$$

We see that the discounted wealth invested in $\theta^0$ is $Y(t) - (Y(t) - \mu)\langle \gamma_t, \mathcal{L}_t p_0 \rangle$ and in $\theta^1$ it is $(Y(t) - \mu)\langle \gamma_t, \mathcal{L}_t p_0 \rangle$. If we want a certain expected return over the period $\mathbb{T}$, then this will of course fix $\mu$ in formula (4.4).

EXAMPLE 4.2. *Exponential utility.* The utility function is

$$(4.8) \qquad u(x) = -\exp(-\mu x),$$

where $\mu > 0$ is given and $x \in \mathbb{R}$. Determination of $\lambda$ gives

$$(4.9) \qquad -\frac{1}{\mu} \ln \frac{\lambda}{\mu} = K_0 + \frac{1}{2\mu} \int_0^{\bar{T}} \sum_{i \in \mathbb{I}} (\Gamma_s^i)^2 \, ds.$$

The optimal discounted wealth process $Y$, for initial wealth $K_0 \in \mathbb{R}$, is given by

$$(4.10) \qquad Y(t) = K_0 - \frac{1}{\mu} \int_0^t \sum_{i \in \mathbb{I}} \Gamma_s^i \, d\bar{W}_s^i.$$

The optimal utility is given by

$$(4.11) \qquad U(t, w) = -\exp\left(-\mu w - \tfrac{1}{2} \int_t^{\bar{T}} \sum_{i \in \mathbb{I}} (\Gamma_s^i)^2 \, ds\right),$$

where $w \in \mathbb{R}$ and $t \in \mathbb{T}$. For an optimal portfolio we get $b_t = -1/\mu$,

$$(4.12) \qquad a_t = (\bar{p}_t(0))^{-1}\left(Y(t) + \frac{1}{\mu}\langle \gamma_t, \mathcal{L}_t p_0\rangle\right)$$

and $\hat{\theta} = \theta^0 + \theta^1$, where $\theta_t^0 = a_t \delta_0$ and

$$(4.13) \qquad \theta_t^1 = -\frac{1}{\mu}\gamma_t(\bar{p}_t)^{-1}\mathcal{L}_t p_0.$$

So in this case the discounted wealth invested in the risky zero-coupon bond of time to maturity $S$ is $\theta_t^1(S)\bar{p}_t(S) = -\frac{1}{\mu}\gamma_t(S)p_0(S+t)$, which is deterministic. However, the portfolio $\theta_t^0$, that is, the number $a_t$ of zero-coupon bonds of time to maturity $0$ is random through its dependence on the discounted wealth $Y(t)$. The discounted wealth invested in $\theta^0$ is $Y(t) + \frac{1}{\mu}\langle \gamma_t, \mathcal{L}_t p_0\rangle$ and in $\theta^1$ it is $-\frac{1}{\mu}\langle \gamma_t, \mathcal{L}_t p_0\rangle$. Expressed in Roll-Overs the portfolio becomes $\hat{\eta} = \eta^0 + \eta^1$,

$$(4.14) \qquad \eta_t^0(S) = Y(t) + \frac{1}{\mu}\langle \gamma_t, \mathcal{L}_t p_0\rangle \delta(S),$$

$$(4.15) \qquad \eta_t^1(S) = -\frac{1}{\mu}\exp\left(\int_0^t (r_s - f_s(S))\,ds\right)\frac{p_0(t+S)}{p_t(S)}\gamma_t(S),$$

where $S \geq 0$ and $\delta = \delta_0$.

EXAMPLE 4.3. *Homogeneous utility.* The utility function is

$$(4.16) \qquad u(x) = x^\mu,$$

where $0 < \mu < 1$ is given and $x > 0$. Determination of $\lambda$ gives

$$(4.17) \qquad \left(\frac{\mu}{\lambda}\right)^{1/(1-\mu)} = K_0 \exp\left(-\frac{\mu}{2(1-\mu)^2}\int_0^{\bar{T}}\sum_{i\in\mathbb{I}}(\Gamma_s^i)^2\,ds\right).$$

The optimal discounted wealth process $Y$, for initial wealth $K_0 > 0$, is given by

$$(4.18) \quad Y(t) = K_0 \exp\left(\int_0^t \sum_{i\in\mathbb{I}}\left(-\frac{1}{1-\mu}\Gamma_s^i\,d\bar{W}_s^i - \frac{1}{2}\left(\frac{1}{1-\mu}\Gamma_s^i\right)^2 ds\right)\right).$$

The optimal utility is given by

$$(4.19) \qquad U(t,w) = w^\mu \exp\left(\frac{\mu}{2(1-\mu)}\int_t^{\bar{T}}\sum_{i\in\mathbb{I}}(\Gamma_s^i)^2\,ds\right), \qquad w > 0.$$

The optimal portfolio $\hat{\theta}$ is given by

$$(4.20) \qquad b_t = -\frac{1}{1-\mu}Y(t)$$

and

$$(4.21) \qquad a_t = (\bar{p}_t(0))^{-1}\Big(1 + \frac{1}{1-\mu}\langle\gamma_t, \mathcal{L}_t p_0\rangle\Big)Y(t),$$

so both $\theta^0$ and $\theta^1$ are proportional to the wealth. The fraction $\theta_t^1(S)\bar{p}_t(S)/Y(t)$ $= -\gamma_t(S)p_0(S+t)/(1-\mu)$, invested in the risky zero-coupon bond of time to maturity $S$, is deterministic.

REMARK 4.4. If, instead of maximizing expected utility of discounted terminal wealth, we maximize expected utility of terminal wealth (see Remark 3.5), we find, in the case of a homogeneous utility function (4.16), that the optimal portfolio $\hat{\varepsilon}$ satisfies

$$\frac{\varepsilon_t^1(S)p_t(S)}{V(t,\hat{\varepsilon})} = \frac{-1}{1-\mu}p_0(t+S)\gamma_t(S) + \frac{-\mu}{1-\mu}\delta_{\bar{T}-t}(S),$$

for time to maturity $S > 0$. The fraction invested in the risky zero-coupon bond of time to maturity $S > 0$ is thus deterministic. In the particular case when the portfolio is restricted to a current account and a zero-coupon bond of maturity exceeding the portfolio management horizon $\bar{T}$, a similar formula was obtained in [10]. It refers to the first term as the Merton result, and to the second as the correction term.

**5. Mathematical complements and proofs.** In the sequel it will be convenient to use a more compact mathematical formalism, which we now introduce. The dual of $\tilde{H}_0$ is identified with $\tilde{H}_0^* = \mathbb{R} \oplus H^*$ by extending the bilinear form, defined in (2.5), to $\tilde{H}_0^* \times \tilde{H}_0$:

$$(5.1) \qquad \langle F, G\rangle = ab + \langle f, g\rangle,$$

where $F = a + f \in \tilde{H}_0^*$, $G = b + g \in \tilde{H}_0$, $a, b \in \mathbb{R}$, $f \in H^*$ and $g \in H$. $\{e_i\}_{i\in\mathbb{N}^*}$ is an orthonormal basis in $\tilde{H}_0$. For $i \in \mathbb{N}^*$, the element $e_i' \in \tilde{H}_0^*$ is given by $\langle e_i', f\rangle = (e_i, f)_{\tilde{H}_0}$, for every $f \in \tilde{H}_0$. The map $[L(E,F)$ denotes the space of linear continuous mappings from $E$ into $F$, $L(E) = L(E,E)]$ $\mathcal{S} \in L(\tilde{H}_0, \tilde{H}_0^*)$ is defined by $\mathcal{S}f = \sum_{i\geq 1}\langle e_i', f\rangle e_i'$. The adjoint $\mathcal{S}^* \in L(\tilde{H}_0^*, \tilde{H}_0)$ is given by $\mathcal{S}^*f = \sum_{i\geq 1}\langle e_i, f\rangle e_i$. Moreover, $(f, g)_{\tilde{H}_0} = \langle \mathcal{S}f, g\rangle$ for $f, g \in \tilde{H}_0$, $(f, g)_{\tilde{H}_0^*} = \langle f, \mathcal{S}^*g\rangle$ for $f, g \in \tilde{H}_0^*$ and $\mathcal{S}$ is unitary. For a given orthonormal basis $\{e_i'\}_{i\in\mathbb{I}}$ in $\tilde{H}_0^*$ we define the $L(\tilde{H}_0)$-valued process $\{\sigma_t\}_{t\in\mathbb{T}}$ by

$$(5.2) \qquad \sigma_t f = \sum_{i\in\mathbb{I}}\langle e_i', f\rangle\sigma_t^i,$$

for $f \in \tilde{H}_0$. We note that if $\sum_{i\geq 1}\|\sigma_t^i\|_{\tilde{H}_0}^2 < \infty$ a.s., then $\sigma_t$ is a.s. a Hilbert–Schmidt operator-valued process, with Schmidt norm

$$(5.3) \qquad \|\sigma_t\|_{H-S} = \Big(\sum_{i\geq 1}\|\sigma_t^i\|_{\tilde{H}_0}^2\Big)^{1/2}.$$

The adjoint is given by

$$\sigma_t^* f = \sum_{i \in \mathbb{I}} \langle f, \sigma_t^i \rangle e_i', \tag{5.4}$$

for $f \in \tilde{H}_0^*$.

We define a cylindrical Wiener process $W$ on $\tilde{H}_0$; cf. Section 4.3.1 of [4]: $W_t = \sum_{i \in \mathbb{I}} W_t^i e_i$. We also define $\Gamma_t = \sum_{i=1}^{\infty} \Gamma_t^i e_i$, which is an element of $\bar{H}_0$ a.s. if $\sum_{i=1}^{\infty} (\Gamma_t^i)^2 < \infty$ a.s. Equation (2.11) now reads

$$\bar{p}_t = \mathcal{L}_t \bar{p}_0 + \int_0^t \mathcal{L}_{t-s} \bar{p}_s m_s \, ds + \int_0^t \mathcal{L}_{t-s} \bar{p}_s \sigma_s \, dW_s, \tag{5.5}$$

its differential

$$d\bar{p}_t = (m_t \bar{p}_t + \partial \bar{p}_t) \, dt + \bar{p}_t \sigma_t \, dW_t, \tag{5.6}$$

equation (2.28)

$$d\bar{G}(t, \theta) = \langle \theta_t, \bar{p}_t m_t \rangle \, dt + \langle \sigma_t^* \bar{p}_t \theta_t, dW_t \rangle, \tag{5.7}$$

relation (2.34)

$$m_t + \sigma_t \Gamma_t = 0 \tag{5.8}$$

and equation (2.36)

$$d\bar{G}(t, \theta) = -\langle \sigma_t^* \bar{p}_t \theta_t, \Gamma_t \rangle \, dt + \langle \sigma_t^* \bar{p}_t \theta_t, dW_t \rangle, \tag{5.9}$$

where $t \in \mathbb{T}$.

The quadratic variation for a process $M$ is, when defined, denoted $\langle\langle M, M \rangle\rangle$.

REMARK 5.1. In order to justify condition (5.8), we note (omitting the a.s.) that if $\theta'$ is a self-financing strategy such that $\theta_t' \in H^*$ is in the annihilator $\{\bar{p}_t \sigma_t^i | i \in \mathbb{I}\}^{\perp} \subset H^*$ of the set $\{\bar{p}_t \sigma_t^i | i \in \mathbb{I}\} \subset H$, then (2.28) gives $d\bar{G}(t, \theta') = \langle \theta_t', m_t \bar{p}_t \rangle \, dt$. $\theta'$ is therefore a riskless self-financing strategy. Since the interest rate of the discounted bank account is zero, in an arbitrage-free market we must have $\langle \theta_t', m_t \bar{p}_t \rangle = 0$. This shows that $m_t \bar{p}_t \in (\{\bar{p}_t \sigma_t^i | i \in \mathbb{I}\}^{\perp})^{\perp}$, that is, $m_t \bar{p}_t$ is an element of the closed linear span of $\{\bar{p}_t \sigma_t^i | i \in \mathbb{I}\}$. Since $\bar{p}_t > 0$, we choose $m_t$ to be an element of the closed linear span $F$ of $\{\sigma_t^i | i \in \mathbb{I}\}$ in $\tilde{H}_0$.

When (also omitting the a.s.) the linear span of $\{\sigma_t^i | i \in \mathbb{I}\}$ has infinite dimension, then condition (5.8) is slightly stronger than $m \in F$, since $\sigma_t$ must be a compact operator in $\tilde{H}_0$. This phenomenon, which is purely due to the infinite dimension of the state space $H$, is not present in the case of a market with a finite number of assets.

REMARK 5.2. The conditions involving $m$ are redundant when equality (5.8) is satisfied. For example, conditions (2.23) and (2.35) imply condition (2.24). In fact, $\|m_t\|_{\tilde{H}_0} \leq (\sum_{i\in\mathbb{I}}\|\sigma_t^i\|_{\tilde{H}_0}^2)^{1/2}(\sum_{i\in\mathbb{I}}|\Gamma_t^i|^2)^{1/2} \leq 1/2(\sum_{i\in\mathbb{I}}\|\sigma_t^i\|_{\tilde{H}_0}^2 + \sum_{i\in\mathbb{I}}|\Gamma_t^i|^2)$. By the Schwarz inequality, $E(\exp(a\int_0^{\bar{T}}\sum_{i\in\mathbb{I}}\|m_t\|_{\tilde{H}_0}\,dt)) \leq (E(\exp(2a\int_0^{\bar{T}}\sum_{i\in\mathbb{I}}\|\sigma_t^i\|_{\tilde{H}_0}^2\,dt)))^{1/2}(E(\exp(2a\int_0^{\bar{T}}\sum_{i\in\mathbb{I}}|\Gamma_t^i|^2\,dt)))^{1/2}$.

REMARK 5.3. When the number of random sources is infinite, that is, $\mathbb{I}=\mathbb{N}^*$, then the straightforward generalization of condition (3.8) from $x\in\mathbb{R}^{\bar{m}}$ to $x\in l^2$ cannot be satisfied, since $\sigma_t$ is a.s. a compact operator in $\tilde{H}_0$. In fact, in this case the left-hand side of (3.8) reads $(x, A(t)x)_{l^2}$. Let $l_t = \mathcal{L}_t p_0$. By the definition of $A(t)$ and by the canonical isomorphism between $l^2$ and $\tilde{H}_0$ we obtain $(x, A(t)x)_{l^2} = \|\sum_{i\in\mathbb{I}}x_i\sigma_t^i l_t\|_{\tilde{H}_0}^2 = \|l_t\sigma_t f\|_{\tilde{H}_0}^2$, where $x_i = (e_i, f)_{\tilde{H}_0}$. Condition (3.8) then reads $\|l_t\sigma_t f\|_{\tilde{H}_0}^2 k_t \geq \|l_t\sigma_t\|_{H-S}\|f\|_{\tilde{H}_0}^2$. Since $\sigma_t$ is a.s. compact, which then also is the case for $l_t\sigma_t$, it follows that $\inf_{\|f\|_{\tilde{H}_0}=1}\|l_t\sigma_t f\|_{\tilde{H}_0} = 0$. This is in contradiction with $k_t$ finite a.s. and $l_t\sigma_t \neq 0$ a.s.

REMARK 5.4. *Concerning condition* (3.11):

(i) $\Gamma$ is unique or more precisely: Given $\sigma$ and $m$ such that the hypotheses of Theorem 3.8 are satisfied, then there is a unique $\Gamma$ satisfying Condition A and satisfying condition (3.11) for some $\gamma$. To establish this fact let $\sigma_t'$ be the usual adjoint operator in $\tilde{H}_0$, with respect to the scalar product in $\tilde{H}_0$, of the operator $\sigma_t$. Condition (3.11) can then be written $\sigma_t'\delta_t = \Gamma_t$, where $\delta_t = \mathcal{S}^* l_t\gamma_t$ and $l_t = \mathcal{L}_t p_0$. $\delta_t \in \tilde{H}_0$, since $\|\delta_t\|_{\tilde{H}_0} = \|l_t\gamma_t\|_{\tilde{H}_0^*} \leq C\|l_t\|_{\tilde{H}_0}\|\gamma_t\|_{\tilde{H}_0^*} < \infty$. This shows that $\Gamma_t$ is in the orthogonal complement, with respect to the scalar product in $\tilde{H}_0$, of $\mathrm{Ker}\,\sigma_t$. There cannot be more than one solution $\Gamma_t$ of (5.8) with this property.

(ii) Condition (3.11) can be satisfied for arbitrary (included degenerated) volatilities $\sigma$, resulting in incomplete markets. An example is obtained by, for given $\sigma$, choosing a $\gamma$ and then defining $\Gamma$ and $m$ by (3.11) and (5.8), respectively.

REMARK 5.5. When $m_t$ and $\sigma_t$ are given functions of $\bar{p}_t$, for every $t\in\mathbb{T}$, then the optimal portfolio problem (3.1) can be considered within a Hamilton–Jacobi–Bellman approach. We illustrate this in the simplest case, when $m_t$ and $\sigma_t$ are deterministic. The optimal value function $U$ then only depends on time $t\in\mathbb{T}$ and on the value on the discounted wealth $w$ at time $t$:

$$(5.10) \qquad U(t,w) = \sup\{E(u(\bar{V}(\bar{T},\theta))|\bar{V}(t,\theta)=w)|\theta\in\mathsf{P}_{\mathsf{sf}}\}$$

[here $E(Y|X=x)$ is the conditional expectation of $Y$ under the condition that $X=x$]. One is then led to the HJB equation

$$
\frac{\partial U}{\partial t}(t,w) + \sup_{f \in H^*} \Bigg\{ -\langle \sigma_t^* f, \Gamma_t \rangle \frac{\partial U}{\partial w}(t,w)
$$
(5.11)
$$
+ \frac{1}{2} \|\sigma_t^* f\|_{H^*}^2 \frac{\partial^2 U}{\partial w^2}(t,w) \Bigg\} = 0
$$

with boundary condition

$$
(5.12) \qquad\qquad U(\bar{T}, w) = u(w).
$$

Equation (5.11) gives

$$
(5.13) \qquad\qquad \frac{\partial U}{\partial t} \frac{\partial^2 U}{\partial w^2} = \frac{1}{2} \|\Gamma_t\|_H^2 \left( \frac{\partial U}{\partial w} \right)^2 .
$$

Each self-financing zero-coupon bond strategy $\hat{\theta} \in \mathsf{P}_{\mathrm{sf}}$, such that

$$
(5.14) \qquad\qquad \langle \hat{\theta}_t, \bar{p}_t \sigma_t^i \rangle = \frac{\Gamma_t^i \, \partial U / \partial w}{\partial^2 U / \partial w^2}, \qquad i \in \mathbb{I},
$$

if $\bar{V}(t, \hat{\theta}) = w$, is then a solution of problem (3.1). In particular, the solution of Corollary 3.9 satisfies (5.14). When $m_t$ and $\sigma_t$ are functions of the price $\bar{p}$, then the HJB equation contains supplementary terms involving Frechét derivatives with respect to $\bar{p}$. The solution of such HJB equations is to our knowledge an open problem.

REMARK 5.6.  *Asymptotic elasticity.* We can prove that there exist utility functions satisfying Condition B with asymptotic elasticity $\limsup_{x \to \infty} x u'(x)/u(x) = 1$. For such $u$, in the situation of Theorems 3.6 and 3.8 there exist optimal portfolios in $\mathsf{P}_{\mathrm{sf}}$. This is in contrast to the situation considered in [11], where for such $u$ and $K_0$ sufficiently large, there is for certain complete financial markets no solution $\hat{X}$ of (3.7) (see Proposition 5.2 of [11]). This remark will be developed in a forthcoming work.

PROOF OF THEOREM 2.1.  Existence, uniqueness and continuity of $\bar{p}$ and $\partial \bar{p}$ follow from Lemma A.1. It then follows from (A.22) of Lemma A.2, with $Y(t) = \mathcal{L}_t \bar{p}_0$, that the solution is given by (2.22). This shows that it is positive.

Finally we prove that condition (2.21) is satisfied. Formula (5.6) and conditions (2.13) and (2.14) give $d(\bar{p}_t(0)) = (\partial \bar{p}_t)(0) \, dt$. Since $\bar{p}_0(0) = 1$, (2.21) follows by integration. □

PROOF OF THEOREM 2.2.  It follows from the explicit expression (2.22) of $\bar{p}$ and (A.20) that $\bar{q} = \tilde{\mathcal{E}}(L)$, where $L(t) = \int_0^t (m_s \, ds + \sum_{i \in \mathbb{I}} \sigma_s^i \, dW_s^i)$. Let $\alpha = 1$ or $\alpha = -1$ and let $J_\alpha = \int_0^t ((\alpha m_s + \alpha(\alpha-1)/2 \sum_{i\in\mathbb{I}} (\sigma_s^i)^2)\,ds + \sum_{i\in\mathbb{I}} \alpha \sigma_s^i\,dW_s^i)$. Then $(\bar{q})^\alpha = \tilde{\mathcal{E}}(J_\alpha)$. According to conditions (2.23) and (2.24), hypotheses (i)–(iv) of Lemma A.4 (with $J_\alpha$ instead of $L$) are satisfied. We now apply estimate (A.40) of Lemma A.4 to $X = (\bar{q})^\alpha$, which proves that $(\bar{q})^\alpha \in L^u(\Omega, P, L^\infty(\mathbb{T}, \tilde{H}_1))$, for $\alpha = \pm 1$. Since $\bar{p}_t = \bar{q}(t)\mathcal{L}_t p_0$, $\mathcal{L}$ is a contraction semigroup and $\tilde{H}_0$ is a Banach algebra, we have $\|\bar{p}_t\|_H^2 + \|\partial\bar{p}_t\|_H^2 \le C(\|\bar{p}_0\|_H^2 + \|\partial\bar{p}_0\|_H^2)\|\bar{q}(t)\|_{\tilde{H}_1}^2$, for some constant $C$ given by $H$. This proves the statement of the lemma in the case $\bar{p}$.

To prove the case of $q$ we note that $\bar{q}(t) = q(t)\bar{p}_t(0)$. Using that the case of $(\bar{q})^\alpha$ is already proved and Hölder's inequality, it is enough to prove that $g \in L^u(\Omega, P, L^\infty(\mathbb{T}, \mathbb{R}))$, where $g(t) = (\bar{p}_t(0))^{-\alpha}$. Since $\bar{p}_t(0) = (\mathcal{L}_t\bar{p}_0)(0)(\bar{q}(t))(0) = \bar{p}_0(t)(\bar{q}(t))(0)$, it follows that $0 \le g(t) = (\bar{p}_0(t))^{-\alpha}((\bar{q}(t))(0))^{-\alpha}$. By Sobolev embedding, $\bar{p}_0$ is a continuous real-valued function on $[0, \infty[$ and it is also strictly positive, so $(\bar{p}_0)^{-\alpha}$ is bounded on $\mathbb{T}$. Once more by Sobolev embedding, $((\bar{q}(t))(0))^{-\alpha} \le C\|(\bar{q}(t))^{-\alpha}\|_{\tilde{H}_0}$. The result now follows, since we have already proved the case of $(\bar{q})^\alpha$. The case of $p$ is so similar to the previous cases that we omit it. $\square$

PROOF OF COROLLARY 2.3. The second part of the proof of Theorem 2.2 gives the result. $\square$

PROOF OF PROPOSITION 2.5. Let $\theta \in \mathsf{P}$ and introduce $X = \sup_{t\in\mathbb{T}} |\bar{G}(t, \theta)|$, $Y(t) = \int_0^t \langle \theta_s, \bar{p}_s m_s \rangle\,ds$ and $Z(t) = \int_0^t \langle \sigma_s^* \bar{p}_s \theta_s, dW_s \rangle$. $\bar{G}(t, \theta) = Y(t) + Z(t)$, according to formula (5.7). Let $\bar{p}$ be given by Theorem 2.1, of which the hypotheses are satisfied.

We shall give estimates for $Y$ and $Z$. By the definition (2.29) of $\mathsf{P}$,

$$(5.15) \qquad E\left(\sup_{t\in\mathbb{T}} (Y(t))^2\right) \le E\left(\left(\int_0^{\bar{T}} |\langle \theta_s, \bar{p}_s m_s \rangle|\,ds\right)^2\right) \le \|\theta\|_{\mathsf{P}}^2.$$

By isometry we obtain

$$
\begin{aligned}
(5.16) \qquad E(Z(t)^2) &= E\left(\int_0^t \left\langle \theta_s \bar{p}_s \sum_{i\in\mathbb{I}} \sigma_s^i\,dW_s^i \right\rangle\right)^2 \\
&= E\left(\int_0^t \sum_{i\in\mathbb{I}} (\langle \theta_s, \bar{p}_s \sigma_s^i \rangle)^2\,ds\right) \\
&\le E\left(\int_0^{\bar{T}} \|\sigma_s^* \theta_s \bar{p}_s\|_{H^*}^2\,ds\right) \le \|\theta\|_{\mathsf{P}}^2.
\end{aligned}
$$

Doob's $L^2$ inequality and inequality (5.16) give $E(\sup_{t\in\mathbb{T}} Z(t)^2) \le 4\|\theta\|_{\mathsf{P}}^2$. Inequality (5.15) then gives $E(X^2) \le 10\|\theta\|_{\mathsf{P}}^2$, which proves the proposition. $\square$

PROOF OF THEOREM 2.8. As we will see, the strong condition (2.35) on $\Gamma$ introduced in (2.34) assures the existence of a martingale measure $Q$ equivalent to $P$, with Radon–Nikodym derivative in $L^u(\Omega, P)$, for each $u \in [1, \infty[$.

LEMMA 5.7. *If* (2.35) *is satisfied, then* $(\xi_t)_{t \in \mathbb{T}}$ *is a* $(P, \mathcal{A})$-*martingale and* $\sup_{t \in \mathbb{T}}(\xi_t)^\alpha \in L^1(\Omega, P)$ *for each* $\alpha \in \mathbb{R}$.

PROOF. Let $M(t) = \int_0^t \sum_{i \in \mathbb{I}} \Gamma_s^i \, dW_s^i$. Then $\langle\langle M, M \rangle\rangle(t) = \int_0^t \sum_{i \in \mathbb{I}} (\Gamma_s^i)^2 \, ds$ and according to condition (2.35) $E(\exp(a \langle\langle M, M \rangle\rangle(\bar{T})) < \infty$, for each $a \geq 0$. By choosing $a = 1/2$, Novikov's criterion (cf. [18], Chapter VIII, Proposition 1.15), shows that $\xi$ is a martingale. Let $b \geq 0$. It then follows from the same reference, by choosing $a = 2b^2$, that $E(\exp(b \sup_{t \in \mathbb{T}} |M(t)|)) < \infty$.

Let $\alpha \in \mathbb{R}$ and let $c(t) = \int_0^t \sum_{i \in \mathbb{I}} (\Gamma_s^i)^2 \, ds$. Then $\xi_t^\alpha = \exp(\alpha M(t) - \alpha/2 c(t))$, so

$$\sup_{t \in \mathbb{T}}(\xi_t)^\alpha \leq \sup_{t \in \mathbb{T}} \exp(|\alpha| M(t) + c(\bar{T})|\alpha|/2)$$

$$\leq \exp\left(\alpha \sup_{t \in \mathbb{T}} |M(t)| + c(\bar{T})|\alpha|/2\right).$$

This and the Schwarz inequality show that

$$\left(E\left(\sup_{t \in \mathbb{T}} \xi_t^\alpha\right)\right)^2 \leq E\left(\exp\left(2|\alpha| \sup_{t \in \mathbb{T}} |M(t)|\right)\right) E(\exp(|\alpha| c(\bar{T}))).$$

The first factor on the right-hand side of this inequality is finite as is seen by choosing $b = 2|\alpha|$, and the second is finite due to condition (2.35). $\square$

The next corollary is a direct application of Girsanov's theorem.

COROLLARY 5.8. *Let* (2.35) *be satisfied. The measure* $Q$, *defined by* $dQ = \xi_{\bar{T}} \, dP$, *is equivalent to* $P$ *on* $\mathcal{F}_{\bar{T}}$ *and* $t \mapsto \bar{W}_t = W_t - \int_0^t \Gamma_s \, ds$, $t \in \mathbb{T}$, *is a cylindrical* $H$-*Wiener process with respect to* $(Q, \mathcal{A})$.

PROOF. According to Lemma 5.7, $\xi$ is a martingale with respect to $(P, \mathcal{A})$. Theorem 10.14 of [4] then gives the result. $\square$

Corollary 5.8 and (2.22) and (5.9) give

$$(5.17) \qquad \bar{p}_t = \exp\left(\int_0^t \mathcal{L}_{t-s}\left(\sum_{i \in \mathbb{I} \sigma_s^i} d\bar{W}_s^i - \tfrac{1}{2} \sum_{i \in \mathbb{I}} (\sigma_s^i)^2 \, ds\right)\right) \mathcal{L}_t \bar{p}_0$$

and

$$(5.18) \qquad d\bar{G}(t, \theta) = \langle \sigma_t^* \bar{p}_t \theta_t, d\bar{W}_t \rangle.$$

To finish the proof of Theorem 2.8, we note that its first part is a restatement of Lemma 5.7 and that its second part is a restatement of Corollary 5.8. □

PROOF OF COROLLARY 2.9. Let $X = \sup_{t \in \mathbb{T}} |\bar{G}(t, \theta)|$. That conditions (2.14) and (2.19) are satisfied follows as in Remark 5.2. The hypotheses of Theorem 2.1 are therefore satisfied and $\bar{p}$ given by Theorem 2.1 is well defined. The square integrability property follows from Proposition 2.5. Finally we have to prove the martingale property. According to hypotheses, (2.35) is satisfied, so Lemma 5.7, Proposition 2.5 and Schwarz inequality give $(E_Q(X))^2 \leq E(\xi_T^2) E(X^2) < \infty$. This shows that $X \in L^1(\Omega, Q)$ and since $\bar{G}(\cdot, \theta)$ is a local $Q$-martingale according to (5.18) it follows that it is a $Q$-martingale (cf. comment after Theorem 4.1 of [18]). □

PROOF OF COROLLARY 2.10. $\bar{V}(t, \theta)$ is given by (2.27), since $\theta \in \mathsf{P}_{\mathrm{sf}}$. According to Corollary 2.9, $\bar{G}(\cdot, \theta)$ is a $Q$-martingale, so this is also the case for $\bar{V}(\cdot, \theta)$. The estimate also follows from Corollary 2.9. We note that if $\bar{V}(\bar{T}, \theta) \geq 0$ and $E_Q(\bar{V}(\bar{T}, \theta)) > 0$, then the martingale property gives $\bar{V}(0, \theta) > 0$, so the market is arbitrage-free. □

PROOF OF LEMMA 3.2. First suppose that $u$ satisfies Condition B. According to condition (3.2) there exists a sufficiently small $x' \in ]\underline{x}, \infty[$, $C > 0$ and $q > 0$ such that for each $x \in ]\underline{x}, x'[$

$$(5.19) \qquad u'(x) \geq C(1 + |x|)^q.$$

With $x = \varphi(y)$ we then have for some $C' > 0$ and for each $y = u'(x) > 0$

$$(5.20) \qquad |\varphi(y)| \leq C' y^{1/q}.$$

Consider case (i). According to condition (3.3) there exist $C > 0$ and $x'' > 0$, such that for each $x \in ]x'', \infty[$

$$(5.21) \qquad u'(x) \leq Cx^{-q}.$$

Then $u'(]\underline{x}, \infty[) = ]0, \infty[$, since $u' > 0$ and according to Condition B(b). With $x = \varphi(y)$, for some $C' > 0$ and for each $y \in ]0, u'(x'')[$

$$(5.22) \qquad |\varphi(y)| \leq C' y^{-1/q}.$$

The continuity of $u'$ and inequalities (5.20) and (5.22) then prove statement (i) with $p = 1/q$.

Consider the case (ii). There exists a unique $x_0 \in ]\underline{x}, \infty[$ such that $u'(x_0) = 0$. Then $u'(x) < 0$ on $]x_0, \infty[$. According to condition (3.4) $\lim_{x \to \infty} u'(x) = -\infty$, so using Condition B(b) we get $u'(]\underline{x}, \infty[) = \mathbb{R}$. Also by (3.4), for some $C > 0$ and $q > 0$, for each $x \in ]x_0, \infty[ \cap ]0, \infty[$, $-u'(x) \geq Cx^q$. We then obtain $0 \le \varphi(y) \le C'|y|^{1/q}$, for some $C'$ and for $y < 0$. This inequality, the continuity of $u'$ and inequality (5.20) then prove statement (ii).

Second, the proof of the converse statement is so similar to the first part of the proof that we omit it. We only note that the definition of $u(\underline{x})$ guarantees that $u$ is u.s.c. □

PROOF OF THEOREM 3.3. We recall that $I = \,]0, \infty[$ if $u' > 0$ on $]\underline{x}, \infty[$ and $I = \mathbb{R}$ if $u'$ takes the value zero in $]\underline{x}, \infty[$, according to Lemma 3.2. We first prove the following lemma:

LEMMA 5.9. *Let $u$ satisfy Condition B and let $\Gamma$ satisfy condition* (2.35). *Then $\varphi(\lambda\xi_{\bar{T}}) \in L^p(\Omega, P)$ for each $p \in [1, \infty[$, $\lambda \in I$, and $\lambda \mapsto E(\xi_{\bar{T}}\varphi(\lambda\xi_{\bar{T}}))$ defines a strictly decreasing homeomorphism from $I$ onto $]\underline{x}, \infty[$. In particular, if $K_0 \in \,]\underline{x}, \infty[$, then there exists a unique $x \in I$ such that $K_0 = E(\xi_{\bar{T}}\varphi(x\xi_{\bar{T}}))$ and $x$ is continuous and strictly decreasing as a function of $K_0$.*

PROOF. Let $\lambda \in I$ and $g_\lambda = \xi_{\bar{T}}\varphi(\lambda\xi_{\bar{T}})$. Lemma 5.7, inequalities (3.5) and (3.6) of Lemma 3.2 and Hölder's inequality show that $\varphi(\lambda\xi_{\bar{T}}) \in L^p(\Omega, P)$ for each $p \in [1, \infty[$. This result and Hölder's inequality give $g_\lambda \in L^1(\Omega, P)$. It follows that $f(\lambda) = E(g_\lambda)$ is well defined.

We show that $f$ is continuous. Let $\{\lambda_n\}_{n \in N^*}$ be a sequence in $I$ converging to $\lambda$. There exists $\bar{\lambda} \in I$ such that $\bar{\lambda} \le \lambda$ and $\bar{\lambda} \le \lambda_n$, for $n \ge 1$. Since $\varphi$ is decreasing and continuous according to Lemma 3.2, we have $|g_{\lambda_n} - g_\lambda| \le 2g_{\bar{\lambda}}$ and $g_{\lambda_n} - g_\lambda \to 0$, a.e. as $n \to \infty$. $g_{\bar{\lambda}} \in L^1(\Omega, P)$, so by Lebesgue's dominated convergence $f(\lambda_n) - f(\lambda) = E(g_{\lambda_n} - g_\lambda) \to 0$, as $n \to \infty$, which proves the continuity.

The function $f$ is decreasing, since $\varphi$ is decreasing. If $\lambda_1, \lambda_2 \in I$ are such that $f(\lambda_1) = f(\lambda_2)$, then $g_{\lambda_1} = g_{\lambda_2}$ a.e. since $\xi_{\bar{T}} > 0$ a.e. $\varphi$ is strictly decreasing, so it follows that $\lambda_1\xi_{\bar{T}} = \lambda_2\xi_{\bar{T}}$. This gives $\lambda_1 = \lambda_2$, which proves that $f$ is strictly decreasing.

The function $\varphi : I \to \,]\underline{x}, \infty[$ is a strictly decreasing bijection, so if $y \to \inf I$ in $I$, then $\varphi(y) \to \infty$ and if $y \to \infty$, then $\varphi(y) \to \underline{x}$. By Fatou's lemma it follows that

$$(5.23) \qquad \liminf_{n \to \infty} f(\lambda_n) \ge E\left(\liminf_{n \to \infty} g_{\lambda_n}\right) = \infty,$$

if $\lambda_n \to \inf I$ in $I$. Let $\lambda_n \to \infty$ in $I$. Choose $\bar{\lambda} \in I$ such that $\bar{\lambda} \le \inf\{\lambda_n | n \ge 1\}$. Then $g_{\bar{\lambda}} - g_{\lambda_n} \ge 0$, since $\varphi$ is decreasing. Application of Fatou's lemma to $g_{\bar{\lambda}} - g_{\lambda_n}$ gives

$$(5.24) \qquad E\left(\limsup_{n \to \infty} g_{\lambda_n}\right) \ge \limsup_{n \to \infty} E(g_{\lambda_n}).$$

If $\underline{x}$ is finite, then (5.24) and, according to Lemma 5.7, $E(\xi_{\bar{T}}) = 1$ give $\underline{x} \geq \limsup_{n \to \infty} E(g_{\lambda_n})$. Since $g_{\lambda_n} \geq \xi_{\bar{T}} \underline{x}$ it follows that

$$\underline{x} = \limsup_{n \to \infty} E(g_{\lambda_n}), \tag{5.25}$$

if $\underline{x}$ is finite. Inequality (5.24) gives

$$-\infty = \limsup_{n \to \infty} E(g_{\lambda_n}) \tag{5.26}$$

if $\underline{x} = -\infty$. Since $f$ is decreasing it follows from (5.23), (5.25) and (5.26) that $f$ is onto $]\underline{x}, \infty[$ and therefore a homeomorphism of $I$ to $]\underline{x}, \infty[$. This completes the proof. $\square$

Now we finish the proof of Theorem 3.3. Let $\mathcal{C}'(K_0) = \{X \in L^2(\Omega, P, \mathcal{F}_{\bar{T}}) | K_0 = E(\xi_{\bar{T}} X)\}$ and let

$$v(x) = \sup_{y \in ]\underline{x}, \infty[} (xy + u(y)), \tag{5.27}$$

$x \in \mathbb{R}$. Here $v$ is the Legendre–Fenchel transform of $-u$. It follows from Condition B that $v \colon \mathbb{R} \to ]-\infty, \infty]$ is l.s.c. and strictly convex; cf. [5]. Let $I = ]0, \infty[$ if $u' > 0$ on $]\underline{x}, \infty[$ and $I = \mathbb{R}$ if $u'$ takes the value zero on $]\underline{x}, \infty[$. Since $-u$ is $C^1$ and strictly convex

$$v(x) = x\varphi(-x) + u(\varphi(-x)), \tag{5.28}$$

for $-x \in I$, which are the elements of the interior of the domain of $v$.

If $\mu \in I$, then $\varphi(\mu\xi_{\bar{T}}) \in L^p(\Omega, P)$ for each $p \in [1, \infty[$, according to Lemma 5.9. Let $\lambda$ be the unique element in $I$, according to Lemma 5.9, such that $\varphi(\lambda\xi_{\bar{T}}) \in \mathcal{C}'(K_0)$. Let $Y = \varphi(\lambda\xi_{\bar{T}})$. Given $X \in \mathcal{C}'(K_0)$. By definition $E(u(X)) = E(u(X)) - \mu(E(\xi_{\bar{T}} X) - K_0)$. It then follows from (5.27) that

$$E(u(X)) = E(u(X) - \lambda\xi_{\bar{T}} X) + \lambda K_0 \leq E(v(-\lambda\xi_{\bar{T}})) + \lambda K_0. \tag{5.29}$$

Formula (5.28) gives that $E(v(-\lambda\xi_{\bar{T}})) = E(u(Y)) - \lambda E(\xi_{\bar{T}} Y)$. Since $Y \in \mathcal{C}'(K_0)$, it follows from (5.29) that

$$E(u(X)) \leq E(u(Y)). \tag{5.30}$$

Therefore $\hat{X} = Y$ is a solution of problem (3.7). This solution is unique since $u$ is strictly concave, which completes the proof. $\square$

PROOF OF COROLLARY 3.4. It follows from Corollary 2.10 that $\{\bar{V}(\bar{T}, \theta) | \theta \in \mathcal{C}(K_0)\} \subset \mathcal{C}'(K_0)$, where $\mathcal{C}'(K_0)$ is given before (5.27). According to Corollary 2.10, $\bar{V}(\cdot, \hat{\theta})$ is a $Q$-martingale, so Theorem 3.3 shows that $\bar{V}(0, \hat{\theta}) = K_0$ and therefore $\hat{\theta} \in \mathcal{C}(K_0)$. This and Theorem 3.3 give

$$\sup_{\theta \in \mathcal{C}(K_0)} E(u(\bar{V}(\bar{T}, \theta))) \leq \sup_{X \in \mathcal{C}'(K_0)} E(u(X)) = E(u(\hat{X})) = E(u(\bar{V}(\bar{T}, \hat{\theta}))),$$

which proves that $\hat{\theta}$ is a solution of problem (3.1). $\square$

PROOF OF THEOREM 3.6. Here $\mathbb{I}$ is a finite set and $K_0 \in ]\underline{x}, \infty[$. We shall construct a portfolio $\hat{\theta} \in \mathcal{C}(K_0)$ such that $\bar{V}(\bar{T}, \hat{\theta}) = \hat{X}$, where $\hat{X}$ is given by Theorem 3.3.

Since $\xi_{\bar{T}}, \hat{X} \in L^p(\Omega, P)$ for each $p \in [1, \infty[$, according to Lemma 5.7 and Theorem 3.3, it follows by Hölder's inequality that $\xi_{\bar{T}} \hat{X} \in L^p(\Omega, P)$, that is, $\hat{X} \in L^p(\Omega, Q)$, for each $p \in [1, \infty[$. In particular $\hat{X} \in L^2(\Omega, Q)$, so by Corollary 5.8 and by the representation of a square integrable random variable as a stochastic integral, there exist progressively measurable real-valued processes $y_i$, $i \in I$, such that $E_Q(\int_0^{\bar{T}} \sum_{i \in \mathbb{I}} y_i(t)^2 \, dt) < \infty$ and such that $\hat{X} = Y(\bar{T})$, where

$$Y(t) = K_0 + \sum_{i \in \mathbb{I}} \int_0^t y_i(s) \, d\bar{W}_s^i, \tag{5.31}$$

for $t \in \mathbb{T}$. We define $y(t) = \sum_{i \in \mathbb{I}} y_i(t) e_i'$. Then $y(t) \in H^*$ a.s. since

$$\|y(t)\|_{H^*}^2 = \sum_{i \in \mathbb{I}} y_i(t)^2. \tag{5.32}$$

Let $Z = \sup_{t \in \mathbb{T}} |Y(t)|$ and let $p \geq 2$. By Doob's inequality, $E_Q(Z^p) \leq (\frac{p}{1-p})^p \times \sup_{t \in \mathbb{T}} E_Q(|Y(t)|^p)$. Now $|Y|^p$ is a $Q$-submartingale, so $E_Q(|Y(t)|^p) \leq E_Q(|Y(\bar{T})|^p) = E_Q(|\hat{X}|^p)$, which gives

$$E_Q(Z^p) < \infty. \tag{5.33}$$

By the Burkholder–Davis–Gundy inequalities, by equality (5.32) and inequality (5.33) one obtains

$$E_Q\left(\left(\int_0^{\bar{T}} \|y(t)\|_{H^*}^2 \, dt\right)^{p/2}\right) < \infty, \tag{5.34}$$

for $p \geq 2$. Since $E(\cdot) = E_Q(\xi_{\bar{T}}^{-1} \cdot)$, it follows from this inequality and from Lemma 5.7 that

$$E\left(\left(\int_0^{\bar{T}} \|y(t)\|_{H^*}^2 \, dt\right)^{p/2}\right) < \infty, \tag{5.35}$$

for $p \geq 2$. We also note that similarly

$$E(Z^p) < \infty, \tag{5.36}$$

for $p \geq 2$.

According to inequality (3.8), $A(t)$ is invertible a.s.; we set $A(t)^{-1} = 0$, when $A(t)$ is not invertible and $A(t)_{ij}^{-1}$ are the matrix elements of $A(t)^{-1}$. We then obtain

$$\|l(t)\sigma_t\|_{H-S} \|A(t)^{-1}\|_{L(\mathbb{R}^{\bar{m}})} \leq Ck(t), \tag{5.37}$$

where $l(t) = \mathcal{L}_t p_0$. Condition (3.8), Schwarz's inequality and inequality (5.37) give

$$(5.38) \qquad E\left(\left(\sup_{t \in \mathbb{T}} \|l(t)\sigma_t\|_{H-S} \|A(t)^{-1}\|_{L(\mathbb{R}^{\bar{m}})}\right)^p\right) < \infty,$$

for $p \in [1, \infty[$.

We define

$$(5.39) \qquad \eta(t) = \sum_{i,j=1}^{\bar{m}} A(t)_{ij}^{-1} l(t) \sigma_t^i y_j(t).$$

It follows from (5.32) that

$$(5.40) \qquad \|\eta(t)\|_H \leq \|A(t)^{-1}\|_{L(\mathbb{R}^{\bar{m}})} \|y(t)\|_{H^*} \|l(t)\sigma_t\|_{H-S},$$

for $t \in \mathbb{T}$. This inequality and inequalities (5.35) and (5.38) give

$$(5.41) \qquad E\left(\left(\int_0^{\bar{T}} \|\eta(t)\|_H^2 \, dt\right)^{p/2}\right) < \infty,$$

for $p \geq 2$. By construction, $\eta(t)$ satisfies

$$(5.42) \qquad (\eta(t), l(t)\sigma_t^i)_H = y_i(t),$$

for $t \in \mathbb{T}$ and $i \in \mathbb{I}$. Defining $\tilde{\theta}_t^1 = \mathcal{S}\eta(t)$, we obtain a solution of the equation

$$(5.43) \qquad \sigma_t^* \tilde{\theta}_t^1 l(t) = y(t),$$

for $t \in \mathbb{T}$. Let $\bar{q}(t) = \bar{p}_t / l(t)$ and $\theta_t^1 = (\bar{q}(t))^{-1} \tilde{\theta}_t^1$. We obtain $\|\theta_t^1\|_{H^*} \leq C \|(\bar{q}(t))^{-1}\|_{\tilde{H}_0} \|\eta(t)\|_H$, where we have used that $\|\tilde{\theta}_t^1\|_{H^*} = \|\eta(t)\|_H$. Theorem 2.2 and inequality (5.41) then give

$$(5.44) \qquad E\left(\left(\int_0^{\bar{T}} \|\theta_t^1\|_{H^*}^2 \, dt\right)^{p/2}\right) < \infty,$$

for $p \geq 2$. Equation (5.43) shows that $\theta^1$ satisfies the equation

$$(5.45) \qquad \sigma_t^* \theta_t^1 \bar{p}_t = y(t),$$

for $t \in \mathbb{T}$. This equality, expression (5.18) of the discounted gains and the martingale representation (5.31), show that

$$(5.46) \qquad Y(t) = K_0 + \bar{G}(t, \theta^1),$$

for $t \in \mathbb{T}$.

We next prove that $\theta^1 \in \mathsf{P}$. By the hypotheses of the theorem it follows that $E((\int_0^{\bar{T}} \sum_{i \in \mathbb{I}} |\Gamma_t^i|^2 \, dt)^{p/2}) < \infty$, for $p \geq 2$. This inequality, definition (2.29) of the portfolio norm, inequality (5.44) and Schwarz inequality give $\|\theta^1\|_{\mathsf{P}} < \infty$, which proves the statement.

Finally we shall construct the announced self-financing strategy $\hat{\theta}$. Let us define the portfolio $\hat{\theta}$ by $\hat{\theta} = \theta^0 + \theta^1$, where $\theta_t^0 = a(t)\delta_0$ and

$$(5.47) \qquad a(t) = ((\bar{p}_t)(0))^{-1}(Y(t) - \langle \theta_t^1, \bar{p}_t \rangle)$$

for $0 \leq t \leq \bar{T}$.

We have to prove that $\hat{\theta} \in \mathcal{C}(K_0)$. To this end it is enough to prove that $\theta^0 \in \mathsf{P}$, since $\theta^1 \in \mathsf{P}$. By definition, we have

$$\|\theta_t^0\|_{H^*} = \sup_{\|f\|_H \leq 1} |\langle \theta_t^0, f \rangle| \leq \sup_{\|f\|_H \leq 1} (|a(t)||f(0)|) \leq C|a(t)|,$$

where the constant is given by Sobolev embedding. Let $b(t) = a(t)\bar{p}_t(0)$. By the definition of $Z$ and Schwarz inequality it follows that

$$\left( E\left( \left( \int_0^{\bar{T}} |b(t)|^2 \, dt \right)^{p/2} \right) \right)^{1/p}$$
$$\leq \bar{T}(E(Z^p))^{1/p} + \left( E\left( \left( \int_0^{\bar{T}} \|\theta_t^1\|_{H^*}^2 \, dt \right)^{p/2} \left( \sup_{t \in \mathbb{T}} \|\bar{p}_t\|_H^p \right) \right) \right)^{1/p},$$

$p \geq 1$. The first term on the right-hand side of this inequality is finite due to (5.36) and the second term is finite due to Theorem 2.2, (5.44) and Schwarz's inequality. Using Corollary 2.3, we obtain now

$$E\left( \left( \int_0^{\bar{T}} |a(t)|^2 \, dt \right)^{p/2} \right) < \infty,$$

$p \geq 1$. This proves in particular that

$$(5.48) \qquad E\left( \int_0^{\bar{T}} \|\theta_t^0\|_{H^*}^2 \, dt \right) < \infty.$$

Since $(\sigma_t)(0) = 0$ according to (2.13), $m_t(0) = 0$ according to (2.14) and by the definition of the norm in $H^*$, we have that $\|\sigma_t^* \theta_t^0 \bar{p}_t\|_{H^*} = 0$ and $\langle \theta_t^0, \bar{p}_t m_t \rangle = 0$. This proves, together with inequality (5.48) and the definition (2.29) of the portfolio norm, that $\theta^0 \in \mathsf{P}$.

We note that by the definition of $\hat{\theta}$, it follows that $\bar{V}(t, \hat{\theta}) = \langle \hat{\theta}(t), \bar{p}_t \rangle = Y(t)$, for $t \in \mathbb{T}$. Moreover, since $(\sigma_t)(0) = 0$, it follows from formula (2.28) that $\bar{G}(t, \theta^0) = 0$, for each $t \in \mathbb{T}$. So by (5.46), $Y(\cdot) = K_0 + \bar{G}(\cdot, \hat{\theta})$, which proves that $\hat{\theta}$ is self-financing with initial value $K_0$.  $\square$

PROOF OF THEOREM 3.8. This proof is, with some exceptions, so similar to the proof of Theorem 3.6 that we only develop the points which are different. Here $\mathbb{I} = \mathbb{N}^*$ or $\mathbb{I} = \{1, \ldots, \bar{m}\}$.

According to Theorem 3.3, there is a unique $\lambda \in I$ such that $\hat{X} = \varphi(\lambda \xi_{\bar{T}})$. Let $M(t) = \int_0^t \sum_{i \in \mathbb{I}} \Gamma_s^i \, dW_s^i$, $t \in \mathbb{T}$. Then $\langle\langle M, M \rangle\rangle$ is deterministic and according to (2.37) and Corollary 5.8, $\xi_t = \exp(M(t) + \frac{1}{2}\langle\langle M, M \rangle\rangle(t))$. Let

$$F(x) = \varphi(\lambda \exp(x + \tfrac{1}{2}\langle\langle M, M \rangle\rangle(\bar{T}))),$$

$x \in \mathbb{R}$. Then $F(M(\bar{T})) = \hat{X}$. We now apply Lemma A.5 to $F$. This gives an integral representation, as in (5.31), with

$$(5.49) \qquad y_i(t) = E_Q(\lambda \xi_{\bar{T}} \varphi'(\lambda \xi_{\bar{T}}) | \mathcal{F}_t) \Gamma_t^i,$$

$i \in \mathbb{I}$ and $t \in \mathbb{T}$.

Using that $\varphi'$ satisfies conditions (3.9) and (3.10) we obtain also here inequalities (5.33) to (5.36).

Let $z(t) = E_Q(\lambda \xi_{\bar{T}} \varphi'(\lambda \xi_{\bar{T}}) | \mathcal{F}_t)$ and let $\gamma$ be given by (3.11). We define $\tilde{\theta}^1 = z\gamma$. By condition (3.11), (5.43) is satisfied.

The remaining part of the proof is the same as for Theorem 3.6. For later reference we observe that $\theta^1 = (l/\bar{p})z\gamma$. $\quad\square$

Proof of Corollary 3.9. The observation in the end of the proof of Theorem 3.8 and expression (5.47) give the stated explicit expression of the optimal portfolio. $\quad\square$

Proof of Theorem 3.10. We first choose a utility function $u$ satisfying Condition C, $u' > 0$ and $\underline{x} = 0$. This is possible as seen by choosing $u(x) = x^{1/2}$, for example. We define $\Theta \in \mathsf{P}_{sf}$ to be the optimal portfolio given by Corollary 3.9 for $K_0 = 1$. Let $\Theta_t = a_t^1 \delta_0 + b_t^1 \gamma_t (\bar{p}_t)^{-1} \mathcal{L}_t p_0$, where $a^1$ and $b^1$ are the coefficients defined by (3.12) and (3.13), respectively. Since $u' > 0$, it follows from Theorem 3.3 and Corollary 3.4 that $\lambda > 0$. It follows from $\lambda \neq 0$, $\varphi' < 0$ and formula (3.12) that $b_t^1 \neq 0$, after a possible redefinition on a set of measure zero.

Since $\underline{x} = 0$, it follows by the definition of $\varphi$ that $\varphi > 0$ and then by Theorem 3.3 and Corollary 3.4 that $\bar{V}(t, \Theta) = E_Q(\bar{V}(\bar{T}, \Theta) | \mathcal{F}_t) > 0$. This shows that statement (i) is satisfied.

Let us now consider a general $u$ satisfying Condition C. The solution $\hat{\theta}$ given by Corollary 3.9, for a general $K_0 \in ]\underline{x}, \infty[$, can now be written

$$\hat{\theta}_t = (a_t - a_t^1 b_t / b_t^1)\delta_0 + (b_t / b_t^1)\Theta_t,$$

which defines $x$ and $y$ in statement (ii) of the theorem. $\quad\square$

## APPENDIX

**A.1. SDEs and $L^p$ estimates.** In this appendix, we state and prove results, used in the article, concerning existence of solutions of some SDEs and $L^p$ estimates of these solutions. Through the Appendix $m$ and $\sigma^i$, $i \in \mathbb{I}$, are $\mathcal{A}$-progressively measurable $\tilde{H}_0$-valued processes satisfying

$$(A.1) \qquad \int_0^{\bar{T}} \left( \|m_t\|_{\tilde{H}_0} + \sum_{i \in \mathbb{I}} \|\sigma_t^i\|_{\tilde{H}_0}^2 \right) dt < \infty \qquad \text{a.s.}$$

The $\tilde{H}_0$-valued semimartingale $L$ is given by

$$(A.2) \qquad L(t) = \int_0^t \left( m_s \, ds + \sum_{i \in \mathbb{I}} \sigma_s^i \, dW_s^i \right) \qquad \text{if } 0 \le t \le \bar{T},$$

and by $L(t) = L(\bar{T})$, if $t > \bar{T}$. We introduce, for $t \ge 0$, the random variable

$$(A.3) \qquad \mu(t) = t + \int_0^t \left( \|m_s\|_{\tilde{H}_0} + \sum_{i \in \mathbb{I}} \|\sigma_s^i\|_{\tilde{H}_0}^2 \right) ds \qquad \text{if } 0 \le t \le \bar{T},$$

and $\mu(t) = t - \bar{T} + \mu(\bar{T})$ if $t > \bar{T}$. $\mu$ is a.s. strictly increasing, absolutely continuous and onto $[0, \infty[$. The inverse $\tau$ of $\mu$ also has these properties and $\tau(t) \le t$. For a continuous $\tilde{H}_0$-valued process $Y$ on $[0, \bar{T}]$ we introduce

$$(A.4) \qquad \rho_t(Y) = \left( E \left( \sup_{s \in [0,t]} \|Y(\tau(s))\|_{\tilde{H}_0}^2 \right) \right)^{1/2},$$

for $t \in [0, \infty[$, where we have defined $Y(t)$ for $t > \bar{T}$ by $Y(t) = Y(\bar{T})$. We note that $\rho_t(Y) \le (E(\sup_{s \in [0,t]} \|Y(s)\|_{\tilde{H}_0}^2))^{1/2}$, since $\tau(t) \le t$.

We will use certain supplementary properties of the Sobolev spaces $H^s$ (cf. Section 7.9 of [8]) and the space $H$. Let $s \ge 0$. There is a norm $\mathcal{N}_s$ equivalent to $\|\cdot\|_{H^s}$, given by

$$(A.5) \quad (\mathcal{N}_s(f))^2 = \sum_{0 \le k \le n} \int_{\mathbb{R}} |f^{(k)}(x)|^2 \, dx + c_s \int_{\mathbb{R}^2} \frac{|f^{(n)}(x) - f^{(n)}(y)|^2}{|x-y|^{1+2s'}} \, dx \, dy,$$

where $f^{(k)}(x) = (d/dx)^k f(x)$, $s = n + s'$, $0 \le s' < 1$, $n \in \mathbb{N}$, $c_{s'} \ge 0$ and $c_0 = 0$. For $f \in H^s$, let $(\kappa f)(x) = f(x)$, $x \ge 0$. The mapping $\kappa \colon H^s \to H^s / H^s_-$ is continuous and surjective, where $H^s_-$ is the closed subspace of $H^s$ of functions with support in $]-\infty, 0]$. Let $\iota \colon H^s / H^s_- \to H^s$ be a continuous linear injective mapping such that $\kappa \iota$ is the identity mapping on $H^s / H^s_-$. To give explicitly such a mapping $\iota$, let $g \in H^s / H^s_-$. For $x \ge 0$, $h$ is defined by $h(x) = g(x)$. For $x < 0$ and $k = 0, \dots, n-1$, let $h^n(x) = (\partial^n g)(-x)$ and $h^k(x) = (\partial^k g)(0) + \int_0^x h^{k+1}(y) \, dy$. Now, for $x < 0$, let $h^k(x) = h^0(x)$. Let $\phi$ be a $C^\infty$ positive function on $\mathbb{R}$, satisfying $\phi(x) = 1$ if $x \ge -1$ and $\phi(x) = 0$ if $x \le -2$. Then $f = h\phi \in H^s$ and $\kappa f = g$. The mapping $\iota g = f$ has the desired properties. In fact it follows using (A.5) and the definition of the norm in $H^s / H^s_-$ that $\mathcal{N}_s(f) \le C_s \|g\|_H$, for some constant $C_s$. Let $\mathbb{R} \ni t \mapsto \mathcal{L}'_t$ be the $C^0$ unitary group of left translations in $H^s$, that is, $(\mathcal{L}'_t f)(x) = f(x+t)$, for $f \in H^s$ and $t, x \in \mathbb{R}$.

Let now $s$ be the given number $s > 1/2$, in (2.6) defining $H$. The map $\kappa$ is extended to $\kappa \colon \mathbb{R} \oplus H^s \to \tilde{H}_0$ by $\kappa(a + f) = a + \kappa f$, where $a \in \mathbb{R}$ and $f \in \tilde{H}_0$. The map $\iota$ is extended to $\iota \colon \tilde{H}_0 \to \mathbb{R} \oplus H^s$ by $\iota(a + f) = a + \iota f$, where $a \in \mathbb{R}$ and $f \in H$. $\mathcal{L}'$ is extended to a $C^0$ unitary group in $\mathbb{R} \oplus H^s$ by $\mathcal{L}'_t(a + f) = a + \mathcal{L}'_t f$, where $t \in \mathbb{R}$, $a \in \mathbb{R}$ and $f \in H$. One easily establishes that with this extended $\mathcal{L}'$

$$(\text{A.6}) \qquad \mathcal{L}_t \kappa = \kappa \mathcal{L}'_t,$$

for all $t \geq 0$.

LEMMA A.1. *If condition* (A.1) *is satisfied and if* $Y$ *is an* $\mathcal{A}$-*progressively measurable* $\tilde{H}_0$-*valued continuous process on* $[0, \bar{T}]$, *satisfying* $\rho_t(Y) < \infty$, *for all* $t \geq 0$, *then the equation*

$$(\text{A.7}) \qquad X(t) = Y(t) + \int_0^t \mathcal{L}_{t-s} X(s) \left( m_s \, ds + \sum_{i \in \mathbb{I}} \sigma_s^i \, dW_s^i \right),$$

$t \in [0, \bar{T}]$, *has a unique solution* $X$, *in the set of* $\mathcal{A}$-*progressively measurable* $\tilde{H}_0$-*valued continuous processes satisfying*

$$(\text{A.8}) \quad \int_0^{\bar{T}} \left( \|X(t)\|_{\tilde{H}_0} + \|X(t) m_t\|_{\tilde{H}_0} + \sum_{i \in \mathbb{I}} \|X(t) \sigma_t^i\|_{\tilde{H}_0}^2 \right) dt < \infty \qquad a.s.$$

*Moreover, this solution satisfies:*

(i) *If* $\int_0^{\bar{T}} (\|m_t\|_{\tilde{H}_1} + \sum_{i \in \mathbb{I}} \|\sigma_t^i\|_{\tilde{H}_1}^2) \, dt < \infty$ *and* $Y$ *is a continuous* $\tilde{H}_1$-*valued process with* $\rho_t(\partial Y) < \infty$, *for all* $t \geq 0$, *then* $X$ *is a continuous* $\tilde{H}_1$-*valued process.*

(ii) *If* (i) *is satisfied and if* $Y$ *is a semimartingale, then* $X$ *is a semimartingale.*

(iii) *If* $Y$ *is* $H$-*valued, then* $X$ *is* $H$-*valued.*

PROOF. The given continuous process $Y$ is extended to $t > \bar{T}$ by $Y(t) = Y(\bar{T})$. For an $\mathcal{A}$-progressively measurable $\tilde{H}_0$-valued process $X$ on $[0, \infty[$, which satisfies (A.8), let

$$(\text{A.9}) \qquad (AX)(t) = \int_0^t \mathcal{L}_{t-s} X(s) \, dL(s),$$

$0 \leq t$. $AX$ is then a continuous process and if $X$ is a solution of (A.7), then $X$ must be a continuous process. It is therefore sufficient to consider existence and uniqueness for continuous process $X$. $A$ is a linear operator from the space of continuous $\tilde{H}_0$-valued processes into itself, since $\tilde{H}_0$ is a Banach algebra, $t \mapsto \mathcal{L}_t$ is a $C^0$ semigroup in $\tilde{H}_0$ and $\mu(\bar{T}) < \infty$ a.s. We note that $(AX)(t)$ is constant for $t \geq \bar{T}$.

It is enough to prove that the equation

$$(\text{A.10}) \qquad X = Y + AX$$

has a unique solution $X$ being a continuous process. Its restriction to $[0, \bar{T}]$ is then the unique solution of the lemma.

In order to introduce the time-transformed equation of (A.10) with respect to $\tau$ let $X'(t) = X(\tau(t))$, $Y'(t) = Y(\tau(t))$, $(A'X')(t) = (AX)(\tau(t))$ and $\rho'_t(X') = (E(\sup_{s \in [0,t]} \|X'(s)\|^2_{\tilde{H}_0}))^{1/2}$. Let also $\mathcal{A}' = (\mathcal{F}_{\tau(t)})_{t \geq 0}$ be the time-transformed filtration. Equation (A.10) has a continuous solution if and only if the time-transformed equation

$$(A.11) \qquad\qquad X' = Y' + A'X'$$

has a continuous solution $X'$.

For given $T > 0$ let $F$ be the Banach space of $\mathcal{A}'$-progressively measurable $\tilde{H}_0$-valued continuous a.s. processes $Z$ on $[0, T]$, with finite norm $\|Z\|_F = \rho'_T(Z)$.

We denote, for $0 \leq t \leq T$, $K_1(t) = \int_0^{\tau(t) \wedge \bar{T}} \mathcal{L}_{\tau(t)-s} X'(\mu(s)) m_s \, ds$ and $K_2(t) = \int_0^{\tau(t) \wedge \bar{T}} \mathcal{L}_{\tau(t)-s} X'(\mu(s)) \sum_{i \in \mathbb{I}} \sigma^i_s \, dW^i_s$, where $a \wedge b = \min\{a, b\}$. Since $\tilde{H}_0$ is an algebra, $\mathcal{L}$ is a $C^0$ contraction semigroup, $\|m_t\|_{\tilde{H}_0} \leq d\mu(t)/dt$ and $\tau$ is the inverse of $\mu$, it follows from Schwarz's inequality that

$$
\begin{aligned}
(\rho'_t(K_1))^2 &\leq CE\left(\left(\int_0^{\tau(t)} \|X'(\mu(s))\|_{\tilde{H}_0} \|m_s\|_{\tilde{H}_0} \, ds\right)^2\right) \\
&\leq CE\left(\left(\int_0^{\tau(t)} \|X'(\mu(s))\|_{\tilde{H}_0} d\mu(s)\right)^2\right) \\
(A.12) \quad &\leq CE\left(\left(\int_0^{\tau(t)} d\mu(s)\right)\left(\int_0^{\tau(t)} \|X'(\mu(s))\|^2_{\tilde{H}_0} d\mu(s)\right)\right) \\
&\leq CtE\left(\int_0^t \|X'(s)\|^2_{\tilde{H}_0} \, ds\right) \\
&\leq CtE\left(\int_0^t \sup_{s' \in [0,s]} \|X'(s')\|^2_{\tilde{H}_0} \, ds\right) \leq Ct \int_0^t (\rho'_s(X'))^2 \, ds,
\end{aligned}
$$

for some $C \geq 0$.

To establish an estimate of $K_2$, we shall use the property (A.6) of the left translation. Since $\kappa$ and $\iota$ are continuous linear operators and $\kappa\iota$ is the identity operator on $\tilde{H}_0$, it follows from (A.6) that

$$(A.13) \qquad K_2(t) = \kappa \mathcal{L}'_{\tau(t)} \int_0^{\tau(t) \wedge \bar{T}} \mathcal{L}'_{-s} \iota X'(\mu(s)) \sum_{i \in \mathbb{I}} \sigma^i_s \, dW^i_s,$$

for all $t \geq 0$. Let $K''_2(t) = \int_0^{\tau(t) \wedge \bar{T}} \mathcal{L}'_{-s} \iota X'(\mu(s)) \sum_{i \in \mathbb{I}} \sigma^i_s \, dW^i_s$, for $t \geq 0$. Then $K''_2$ is an $\mathbb{R} \oplus H^s$-valued square integrable martingale, with respect to the time-transformed filtration $\mathcal{A}'$. In fact, we obtain by isometry, the unitarity of $\mathcal{L}'$ and as in the case of $K_1$, that

$$
\begin{aligned}
E(\|K_2''(t)\|_{\mathbb{R}\oplus H^s}^2) &\le E\left(\int_0^{\tau(t)} \left\| \mathcal{L}'_{-u}\iota X(u) \sum_{i\in\mathbb{I}} \sigma_u^i \right\|_{\mathbb{R}\oplus H^s}^2 du\right) \\
&\le CE\left(\int_0^{\tau(t)} \|X'(\mu(u))\|_{\tilde{H}_0}^2 \sum_{i\in\mathbb{I}} \|\sigma_u^i\|_{\tilde{H}_0}^2 \, du\right) \\
&\le CE\left(\int_0^t \sup_{u'\in[0,u]} \|X'(\mu(u'))\|_{\tilde{H}_0}^2 \, d\mu(u)\right) \\
&\le C\int_0^t (\rho_u'(X'))^2 \, du,
\end{aligned}
$$

(A.14)

for some $C > 0$ and for all $t \ge 0$. Since $\mathcal{L}'_t$ is unitary and a $\kappa$ is continuous with norm equal to 1, it follows from (A.13) that $(\rho_t'(K_2))^2 \le E(\sup_{u\in[0,t]} \|K_2''(u)\|_{\mathbb{R}\oplus H^s}^2)$. By Doob's inequality (cf. Theorem 3.8 of [4]) we have $E(\sup_{u\in[0,t]} \|K_2''(u)\|_{\mathbb{R}\oplus H^s}^2) \le 4\sup_{u\in[0,t]} E(\|K_2''(u)\|_{\mathbb{R}\oplus H^s}^2)$. This gives, together with inequality (A.14), that

$$
(A.15) \qquad (\rho_t'(K_2))^2 \le C\int_0^t (\rho_s'(X'))^2 \, ds,
$$

for $t \ge 0$, where $C$ chosen sufficiently big is independent of $t$. Formula (A.9) and inequalities (A.12) and (A.15) show that for $t \in [0, T]$,

$$
(A.16) \qquad (\rho_t'(A'X'))^2 \le C'^2(1+t)\int_0^t (\rho_s'(X'))^2 \, ds,
$$

where $C'$ is a constant independent of $T$. In particular,

$$
(A.17) \qquad \rho_t'(A'X') \le C'(1+t)^{1/2}t^{1/2}\rho_t'(X'),
$$

$t \in [0, T]$.

If $T > 0$ is sufficiently small, then (A.17) gives that $\|A'X'\|_F \le a\|X'\|_F$, where $0 \le a < 1$. Therefore $A' \in L(F)$ and $I + A'$ has bounded inverse. Equation (A.11) has then a unique solution $X' \in F$. Let $h(t) = \int_0^t (\rho_s'(X'))^2 \, ds$ and $a(t) = \int_0^t (\rho_s'(Y'))^2 \, ds$. Equation (A.11) and inequality (A.16) show that a solution $X' \in F$ satisfies

$$
(A.18) \qquad h(t) \le 2a(t) + 2C'^2\int_0^t (1+s)h(s) \, ds,
$$

for $t \in [0, T]$. Grönwall's inequality gives $h(t) \le 2a(t)\exp(C'^2 t(2+t))$. Equation (A.11) and inequality (A.16) then show that there exists a finite constant $C_T''$ for every $T > 0$ independent of $X'$, such that $\|X'\|_F^2 = (\rho_T'(X'))^2 \le C_T''(\rho_T'(Y'))^2$. It follows that the solution can be extended to all $T > 0$ and this extended solution is unique. This proves the statement of the lemma concerning the existence and uniqueness of an $\tilde{H}_0$-valued continuous solution of (A.7).

We next prove the supplementary statements (i), (ii) and (iii):

(i) We have just to replace, in the above proof, the space $\tilde{H}_0$ by $\tilde{H}_1$ and redefine appropriately the maps $\iota$ and $\kappa$.

(ii) If $Y$ is a semimartingale, then Itô's lemma and the fact that $\partial X$ is a continuous process give

$$(A.19) \qquad dX(t) = dY(t) + \partial(X(t) - Y(t))\, dt + X(t)\, dL(t).$$

This shows that $X$ is a semimartingale.

(iii) $H$ is a closed subspace of $\tilde{H}_0$ and if $X$ is $H$-valued, then $AX$ is also $H$-valued. This shows that the unique solution of (A.10) is $H$-valued. $\quad\square$

The solution of (A.7) can be given explicitly, which we shall use to derive estimates of the solution. Let

$$(A.20) \quad (\tilde{\mathcal{E}}(L))(t) = \exp\left( \int_0^t \mathcal{L}_{t-s}\left( \left( m_s - \tfrac{1}{2} \sum_{i \in \mathbb{I}} (\sigma_s^i)^2 \right) ds + \sum_{i \in \mathbb{I}} \sigma_s^i\, dW_s^i \right) \right),$$

for $t \in \mathbb{T}$.

LEMMA A.2. *Let condition* (A.1) *be satisfied. Then* $\tilde{\mathcal{E}}(L)$ *is the unique* $\tilde{H}_0$*-valued continuous a.s. solution of*

$$(A.21) \qquad (\tilde{\mathcal{E}}(L))(t) = 1 + \int_0^t \mathcal{L}_{t-s}(\tilde{\mathcal{E}}(L))(s)\, dL(s),$$

*for* $t \in \mathbb{T}$. *Let also* $L'(t) = \int_0^t \sum_{i \in \mathbb{I}} (\sigma_s^i)^2\, ds - L(t)$. *Then the unique solution* $X$ *of* (A.7) *in Lemma* A.1 *is given by*

$$(A.22) \qquad X(t) = Y(t) - (\tilde{\mathcal{E}}(L))(t) \int_0^t \mathcal{L}_{t-s} Y(s)(\tilde{\mathcal{E}}(L'))(s)\, dL'(s),$$

*for* $t \in \mathbb{T}$.

PROOF. Let $l_T(t) = \int_0^t \mathcal{L}_{T-s}(m_s + \sum_{i \in \mathbb{I}} \sigma_s^i\, dW_s^i)$, let $n_T(t) = \int_0^t \mathcal{L}_{T-s}((m_s - \tfrac{1}{2} \sum_{i \in \mathbb{I}} (\sigma_s^i)^2)\, ds + \sum_{i \in \mathbb{I}} \sigma_s^i\, dW_s^i)$, for $0 \le t \le T \le \tilde{T}$ and let $N(t) = n_t(t)$, for $t \in \mathbb{T}$. Then $N$ is an $\tilde{H}_0$-valued continuous process, according to the hypothesis on $m$ and $\sigma$ and since $\tilde{H}_0$ is a Banach algebra. This is then also the case of $\tilde{\mathcal{E}}(L)$, since $\|(\tilde{\mathcal{E}}(L))(t)\|_{\tilde{H}_0} \le \exp(\|N(t)\|_{\tilde{H}_0})$. We note that $dl_T(t) = \mathcal{L}_{T-t}\, dL(t)$ and that $\mathcal{L}_{T-t}(\tilde{\mathcal{E}}(L))(t) = \exp(n_T(t))$. Integration gives

$$(A.23) \qquad \begin{aligned} \int_0^t \mathcal{L}_{t-s}(\tilde{\mathcal{E}}(L))(s)\, dL(s) &= \int_0^t \exp(n_t(s))\, dl_t(s) \\ &= \exp(n_t(t)) - 1 = \tilde{\mathcal{E}}(L)(t) - 1. \end{aligned}$$

This proves that $\tilde{\mathcal{E}}(L)$ is a solution of (A.21). The uniqueness follows from Lemma A.1.

To prove (A.22), let $l'_T(t) = \int_0^t \mathcal{L}_{T-s} \sum_{i \in \mathbb{I}} (\sigma_s^i)^2 \, ds - l_T(t)$, for $0 \le t \le T \le \bar{T}$. Then

(A.24)
$$d \exp(n_T(t)) = \exp(n_T(t)) \, dl_T(t) \quad \text{and}$$
$$d \exp(-n_T(t)) = \exp(-n_T(t)) \, dl'_T(t).$$

Let also $y_T(t) = \mathcal{L}_{T-s} Y(t)$ and $z_T(t) = \mathcal{L}_{T-s}(X(t) - Y(t))/(\tilde{\mathcal{E}}(L))(t)$, for $0 \le t \le T \le \bar{T}$. Let $X$ be the unique solution given by Lemma A.1. Equation (A.7) then reads $X(t) = Y(t) + \int_0^t \mathcal{L}_{t-s} X(s) \, dL(s)$, for $t \in \mathbb{T}$. Applying $\mathcal{L}_{T-t}$, with $0 \le t \le T \le \bar{T}$, on both sides we obtain

(A.25) $\quad z_T(t) = \exp(-n_T(t)) \int_0^t (y_T(s) + z_T(s) \exp(n_T(s))) \, dl_T(s).$

Itô's lemma and formulas (A.24), (A.25), $z_T(0) = 0$ give

$$z_T(t) = \int_0^t \exp(-n_T(s))(y_T(s) + z_T(s) \exp(n_T(s))) \, dl_T(s)$$
$$+ \int_0^t z_T(s) \, dl'_T(s) - \int_0^t \mathcal{L}_{T-s} \sum_{i \in \mathbb{I}} (\sigma_s^i)^2 (\exp(-n_T(s)) y_T(s) + z_T(s)) \, ds,$$

for $0 \le t \le T \le \bar{T}$. Rewriting this formula we obtain

$$z_T(t) = \int_0^t z_T(s) \left( dl_T(s) + dl'_T(s) - \mathcal{L}_{T-s} \sum_{i \in \mathbb{I}} (\sigma_s^i)^2 \, ds \right)$$
$$+ \int_0^t y_T(s) \exp(-n_T(s)) \left( dl_T(s) - \mathcal{L}_{T-s} \sum_{i \in \mathbb{I}} (\sigma_s^i)^2 \, ds \right).$$

The definitions of $l$ and $l'$ then give

(A.26) $\qquad z_T(t) = - \int_0^t y_T(s) \exp(-n_T(s)) \, d'l_T(s),$

for $0 \le t \le T \le \bar{T}$. Choosing $T = t$, we now obtain equation (A.22) since $y_T(s) \exp(-n_T(s)) \, d'l_T(s) = \mathcal{L}_{T-s} Y(s)(\tilde{\mathcal{E}}(L'))(s) \, dL'(s)$, for $0 \le s \le T \le \bar{T}$. $\qquad \square$

The next technical lemma collects estimates of norms of certain $\tilde{H}_0$-valued processes that we need later.

LEMMA A.3. *Let*

$$\|(m, \sigma)\|_j = \int_0^{\bar{T}} \sum_{0 \le k \le j} \|\partial^k m_t\|_{\tilde{H}_0} \, dt + \left( \int_0^{\bar{T}} \sum_{0 \le k \le j} \|\partial^k \sigma_t\|_{H-S}^2 \, dt \right)^{1/2},$$

$j \in \mathbb{N}$, and let $Z(t) = \int_0^t \mathcal{L}_{t-s} \, dL(s)$, $t \in \mathbb{T}$. Let $F : [0, \infty[ \to [0, \infty[$ be a function which is continuous together with its first two derivatives and which has $F' \geq 0$.

(i) If $\|(m, \sigma)\|_0 < \infty$, then

$$(A.27) \qquad F(\|Z(t)\|_{\tilde{H}_0}^2) \leq F(0) + \int_0^t \left( a(s) \, ds + \sum_{i \in \mathbb{I}} b_i(s) \, dW_s^i \right),$$

where $a$ and $b_i$, $i \in \mathbb{I}$, are progressively measurable processes satisfying

$$(A.28) \qquad \begin{aligned} |a(t)| &\leq F'(\|Z(t)\|_{\tilde{H}_0}^2)(2\|Z(t)\|_{\tilde{H}_0} \|m_t\|_{\tilde{H}_0} + \|\sigma_t\|_{H-S}^2) \\ &\quad + 2|F''(\|Z(t)\|_{\tilde{H}_0}^2)| \|Z(t)\|_{\tilde{H}_0}^2 \|\sigma_t\|_{H-S}^2 \end{aligned}$$

and

$$(A.29) \qquad b_i(t) = 2F'(\|Z(t)\|_{\tilde{H}_0}^2)(Z(t), \sigma_t^i)_{\tilde{H}_0},$$

$t \in \mathbb{T}$.

(ii) Moreover, if $\|(m, \sigma)\|_1 < \infty$, then

$$(A.30) \qquad \begin{aligned} F(\|Z(t)\|_{\tilde{H}_0}^2) = F(0) + \int_0^t \bigg( &(-v(s)F'(\|Z(s)\|_{\tilde{H}_0}^2) + a(s)) \, ds \\ &+ \sum_{i \in \mathbb{I}} b_i(s) \, dW_s^i \bigg), \end{aligned}$$

where

$$(A.31) \qquad v(t) = -2(Z(s), \partial Z(s))_{\tilde{H}_0} \geq 0,$$

with $t \in \mathbb{T}$.

PROOF. Suppose first that $\|(m, \sigma)\|_j < \infty$, for each $j \in \mathbb{N}$. We remember that the set $\mathcal{D}_\infty$ of all $f \in \tilde{H}_0$, such that $\partial^j f \in \tilde{H}_0$ for each $j \in \mathbb{N}$, is dense in $\tilde{H}_0$. Then $Z(t) \in \mathcal{D}_\infty$ a.s. Itô's lemma gives

$$(A.32) \quad \|Z(t)\|_{\tilde{H}_0}^2 = \int_0^t \left( (2(Z(s), \partial Z(s))_{\tilde{H}_0} + a^{(1)}(s)) \, ds + \sum_{i \in \mathbb{I}} b_i^{(1)}(s) \, dW_s^i \right),$$

where

$$(A.33) \qquad a^{(1)}(t) = 2(Z(t), m_t)_{\tilde{H}_0} + \|\sigma_t\|_{H-S}^2$$

and

$$(A.34) \qquad b_i^{(1)}(t) = 2(Z(t), \sigma_t^i)_{\tilde{H}_0}.$$

We note that $|a^{(1)}(t)| \leq 2\|Z(t)\|_{\tilde{H}_0}\|m_t\|_{\tilde{H}_0} + \|\sigma_t\|_{\tilde{H}-S}^2$ and that $|b_i^{(1)}(t)| \leq 2\|Z(t)\|_{\tilde{H}_0}\|\sigma_t^i\|_{\tilde{H}_0}$. Once more, by Itô's lemma we obtain

$$
\begin{aligned}
F(\|Z(t)\|_{\tilde{H}_0}^2) = F(0) \\
(A.35) \qquad + \int_0^t \Bigg( \big(2(Z(s), \partial Z(s))_{\tilde{H}_0} F'(\|Z(s)\|_{\tilde{H}_0}^2) + a(s)\big) ds \\
+ \sum_{i \in \mathbb{I}} b_i(s)\, dW_s^i \Bigg),
\end{aligned}
$$

where

$$
(A.36) \qquad a(t) = F'(\|Z(t)\|_{\tilde{H}_0}^2) a^{(1)}(t) + \tfrac{1}{2} F''(\|Z(t)\|_{\tilde{H}_0}^2) \sum_{i \in \mathbb{I}} (b_i^{(1)}(t))^2
$$

and

$$
(A.37) \qquad b_i(t) = F'(\|Z(t)\|_{\tilde{H}_0}^2) b_i^{(1)}(t).
$$

Inequality (A.28) of the lemma follows from the noted estimates for $a^{(1)}(t)$ and $b_i^{(1)}(t)$, from (A.36) and from $F' \geq 0$. Formula (A.29) of the lemma follows from (A.34) and (A.37). Inequality (A.27) of the lemma follows from equality (A.30) of the lemma. Equality (A.30) follows from equality (A.35) and the definition of $v$ in (A.31). $Z(t) \in \mathcal{D}_\infty$ and $\partial$ is the generator of a $C^0$ contraction semigroup, in a real Hilbert space, which give the inequality in (A.31).

We have now proved all statements of the lemma under the supplementary hypothesis that $\|(m,\sigma)\|_j < \infty$, for each $j \in \mathbb{N}$. The general case is now obtained by continuity. $\quad\square$

In the next lemma we establish that the solution of (A.7) is in $L^p$, $p \in [0, \infty[$.

LEMMA A.4.   *Let condition* (A.1) *be satisfied and let* (i) $E(\exp(p \int_0^{\bar{T}} (\|m_t\|_{\tilde{H}_0} + \sum_{i \in \mathbb{I}} \|\sigma_t^i\|_{\tilde{H}_0}^2)\, dt)) < \infty$, *for each* $p \in [1, \infty[$. *Suppose that* $Y$ *in Lemma* A.1 *satisfies* (ii) $E(\sup_{t \in \mathbb{T}} \|Y(t)\|_{\tilde{H}_0}^p) < \infty$, *for each* $p \in [1, \infty[$. *Then the unique solution* $X$ *of* (A.7) *in Lemma* A.1 *satisfies*

$$
(A.38) \qquad E\left(\sup_{t \in \mathbb{T}} \|X(t)\|_{\tilde{H}_0}^p\right) < \infty,
$$

*for each* $p \in [1, \infty[$. *In particular, if* $\tilde{\mathcal{E}}(L)$ *is as in Lemma* A.2, *then*

$$
(A.39) \qquad E\left(\sup_{t \in \mathbb{T}} \|\tilde{\mathcal{E}}(L)\|_{\tilde{H}_0}^p\right) < \infty,
$$

*for each* $p \in [1, \infty[$. *Moreover, if* (iii) $E((\int_0^{\bar{T}}(\|m_t\|_{\tilde{H}_1} + \sum_{i \in \mathbb{I}} \|\sigma_t^i\|_{\tilde{H}_1}^2) dt)^p) < \infty$ *and* (iv) $E(\sup_{t \in \mathbb{T}} \|Y(t)\|_{\tilde{H}_1}^p) < \infty$, *for each* $p \in [1, \infty[$, *then also*

$$(A.40) \qquad E\left(\sup_{t \in \mathbb{T}} \|X(t)\|_{\tilde{H}_1}^p\right) < \infty,$$

*for each* $p \in [1, \infty[$.

Proof. Suppose that conditions (i) and (ii) are satisfied.

We first prove inequality (A.39). Let $N(t) = \int_0^t \mathcal{L}_{t-s}((m_s - \frac{1}{2} \sum_{i \in \mathbb{I}} (\sigma_s^i)^2) ds + \sum_{i \in \mathbb{I}} \sigma_s^i dW_s^i)$, for $t \in \mathbb{T}$. Then $\tilde{\mathcal{E}}(L) = \exp(N(t))$ according to (A.20). Since $\tilde{H}_0$ is a Banach algebra it follows that

$$(A.41) \quad \|(\tilde{\mathcal{E}}(L))(t)\|_{\tilde{H}_0} \leq \exp(C\|N(t)\|_{\tilde{H}_0}) \leq \exp(C(1 + \|N(t)\|_{\tilde{H}_0}^2)^{1/2}),$$

for a constant $C$ given by $\tilde{H}_0$.

We use Lemma A.3 to find a bound of the right-hand side of (A.41). Let $a$ and $b_i$ be given by Lemma A.3, with $F(x) = (1 + x)^{1/2}$, let $A(t) = \int_0^t |a(s)| ds$ and let $M(t) = \int_0^t \sum_{i \in \mathbb{I}} b_i(s) dW_s^i$. Then inequality (A.27) gives

$$(A.42) \qquad (1 + \|N(t)\|_{\tilde{H}_0}^2)^{1/2} \leq 1 + A(t) + M(t),$$

inequality (A.28) gives

$$(A.43) \qquad |a(t)| \leq \|m_t\|_{\tilde{H}_0} + \frac{3}{2}C\|\sigma_t\|_{H-S}^2$$

and (A.29) gives

$$(A.44) \qquad b_i(t) = (1 + \|N(t)\|_{\tilde{H}_0}^2)^{-1/2}(N(t), \sigma_t^i)_{\tilde{H}_0},$$

where $i \in \mathbb{I}$ and $t \in \mathbb{T}$. Obviously $|b_i(t)| \leq \|\sigma_t^i\|_{\tilde{H}_0}$ and the quadratic variation $\langle\langle M, M \rangle\rangle(t) \leq \int_0^t \|\sigma_s\|_{H-S}^2 ds$.

By the hypothesis of the lemma and (A.43) it follows that

$$(A.45) \qquad E(\exp(pA(\bar{T}) + p\langle\langle M, M \rangle\rangle(\bar{T}))) < \infty$$

for each $p \in [1, \infty[$. Novikov's criteria (cf. [18], Chapter VIII, Proposition 1.15) and inequality (A.45) give

$$(A.46) \qquad E\left(\exp\left(p \sup_{t \in \mathbb{T}} |M(t)|\right)\right) < \infty$$

for each $p \in [1, \infty[$. Inequality (A.42) gives

$$(A.47) \quad E(\exp(q(1 + \|N(t)\|_{\tilde{H}_0}^2)^{1/2})) \leq E\left(\exp\left(q\left(1 + A(\bar{T}) + \sup_{t \in \mathbb{T}} |M(t)|\right)\right)\right)$$

for each $q \in [0, \infty[$. It follows from Schwarz's inequality and inequalities (A.45), (A.46) and (A.47) that

$$(A.48) \qquad E(\exp(q(1 + \|N(t)\|_{\tilde{H}_0}^2)^{1/2})) < \infty$$

for each $q \in [0, \infty[$. Statement (A.39) now follows from inequalities (A.41) and (A.48) by choosing $q = pC$.

We use the explicit expression (A.22) for $X$ to prove (A.38). Let $Z(t) = \int_0^t \mathcal{L}_{t-s} \, dV(s)$, where $V(t) = \int_0^t Y(s)(\tilde{\mathcal{E}}(L'))(s) \, dL'(s)$ and $L'$ is as in Lemma A.2. Explicitly

$$V(t) = \int_0^t \left( \alpha(s) \, ds + \sum_{i \in \mathbb{I}} \beta_i(s) \, dW_s^i \right),$$

where $\alpha(t) = Y(t)(\tilde{\mathcal{E}}(L'))(t)((\sum_{i \in \mathbb{I}} \sigma_t^i)^2 - m_t)$ and $\beta_i(t) = -Y(t)(\tilde{\mathcal{E}}(L'))(t)\sigma_t^i$. Since we have proved (A.39), by Schwarz's inequality it is enough to prove

$$(A.49) \qquad E\left( \sup_{t \in \mathbb{T}} \|Z(t)\|_{\tilde{H}_0}^p \right) < \infty,$$

for each $p \in [1, \infty[$, to establish (A.38). We proceed as we did earlier in this proof to obtain (A.42). We now obtain using Lemma A.3

$$(A.50) \qquad (1 + \|Z(t)\|_{\tilde{H}_0}^2)^{1/2} \le 1 + A_1(t) + M_1(t),$$

where $A_1(t) = \int_0^t |a_1(s)| \, ds$, $M_1(t) = \int_0^t \sum_{i \in \mathbb{I}} b_{1i}(s) \, dW_s^i$,

$$(A.51) \quad \begin{aligned} |a_1(t)| &\le C' \|Y(t)\|_{\tilde{H}_0} \|(\tilde{\mathcal{E}}(L'))(t)\|_{\tilde{H}_0} \\ &\quad \times (\|m_t\|_{\tilde{H}_0} + (1 + \|Y(t)\|_{\tilde{H}_0} \|(\tilde{\mathcal{E}}(L'))(t)\|_{\tilde{H}_0}) \|\sigma_t\|_{H-S}^2), \end{aligned}$$

with $C'$ given by $H$ and

$$(A.52) \qquad b_{1i}(t) = -(1 + \|Z(t)\|_{\tilde{H}_0}^2)^{-1/2} (Z(t), Y(t)(\tilde{\mathcal{E}}(L'))(t)\sigma_t^i)_{\tilde{H}_0}.$$

Choosing the constant $C'$ sufficiently big, (A.52) gives

$$(A.53) \quad \langle\langle M_1, M_1 \rangle\rangle(t) \le C' \int_0^t \|Y(s)\|_{\tilde{H}_0}^2 \|(\tilde{\mathcal{E}}(L'))(s)\|_{\tilde{H}_0}^2 \|\sigma_s\|_{H-S}^2 \, ds.$$

Hölder's inequality, inequalities (A.42) and (A.51) and the hypotheses of the lemma give

$$(A.54) \qquad E\left( \sup_{t \in \mathbb{T}} (A_1(t))^p \right) < \infty,$$

for each $p \in [1, \infty[$. Similarly, using $\int_0^t \|Y(s)\|_{\tilde{H}_0}^2 \|(\tilde{\mathcal{E}}(L'))(s)\|_{\tilde{H}_0}^2 \|\sigma_s\|_{H-S}^2 \, ds$ $\le (\sup_{s \in \mathbb{T}} \|Y(s)\|_{\tilde{H}_0}^2)(\sup_{s \in \mathbb{T}} \|(\tilde{\mathcal{E}}(L'))(s)\|_{\tilde{H}_0}^2) \int_0^t \|\sigma_s\|_{H-S}^2 \, ds$, (A.53) gives

$$(A.55) \qquad E((\langle\langle M_1, M_1 \rangle\rangle(\bar{T}))^{p/2}) < \infty,$$

for each $p \in [1, \infty[$. The BDG inequality then gives

$$(A.56) \qquad E\left(\sup_{t \in \mathbb{T}} (M_1(t))^p\right) < \infty,$$

for each $p \in [1, \infty[$. Now inequalities (A.50), (A.54) and (A.56) prove (A.49).

Finally, to prove inequality (A.40), we suppose also that conditions (iii) and (iv) are satisfied.

The solution $X$ of (A.7) is, according to Lemma A.1, in the domain of $\partial$, that is, $X(t) \in \tilde{H}_1$. Since $\partial$ is continuous from $\tilde{H}_1$ to $\tilde{H}_0$, we have $\partial \int_0^t \mathcal{L}_{t-s} X(s) \, dL(s) = \int_0^t \mathcal{L}_{t-s} \partial X(s) \, dL(s)$. Application of $\partial$ on both sides of (A.7) then gives

$$(A.57) \qquad X_1(t) = Y_1(t) + \int_0^t \mathcal{L}_{t-s} X_1(s) \, dL(s),$$

where $X_1(t) = \partial X(t)$, $Y_1(t) = \partial Y(t) + \int_0^t \mathcal{L}_{t-s} X(s) \, dL_1(s)$, with $L_1(t) = \int_0^t (\partial m_s \, ds + \sum_{i \in \mathbb{I}} \partial \sigma_s^i \, dW_s^i)$. We can now use inequality (A.38) for $X_1$, since in the context of (A.57) hypotheses (i) and (ii) are satisfied. This proves inequality (A.40). $\square$

For completeness we prove, for the case of an infinite number of random sources, a representation result. The measure $Q$ and the cylindrical Wiener process $\bar{W}$ are as in Corollary 5.8.

LEMMA A.5. *Let $\Gamma$ be deterministic and satisfy condition*

$$(A.58) \qquad \int_0^{\bar{T}} \sum_{i \in \mathbb{I}} |\Gamma_t^i|^2 \, dt < \infty \qquad a.s.$$

*and let $M(t) = \int_0^t \sum_{i \in \mathbb{I}} \Gamma_s^i \, d\bar{W}_s^i$, $t \in \mathbb{T}$. If $F \in C(\mathbb{R})$ is absolutely continuous, with derivative $F'$, and $E_Q(F(M(\bar{T}))^2 + F'(M(\bar{T}))^2) < \infty$, then*

$$(A.59) \quad F(M(\bar{T})) = E_Q(F(M(\bar{T}))) + \int_0^{\bar{T}} E_Q(F'(M(\bar{T})) | \mathcal{F}_t) \, dM(t),$$

*for each $t \in \mathbb{T}$.*

PROOF. We have $\langle\langle M, M \rangle\rangle(t) = \int_0^t \sum_{i \in \mathbb{I}} (\Gamma_s^i)^2 \, ds < \infty$ according to condition (A.58) and the quadratic variation $\langle\langle M, M \rangle\rangle$ is deterministic. Let $n_{\mu,t}(x) = \exp(ix\mu + \frac{\mu^2}{2} \langle\langle M, M \rangle\rangle(t))$, $\mu \in \mathbb{R}$, and let $n'_{\mu,t}(x)$ be the derivative with respect to $x$ of $n_{\mu,t}(x)$. Then $\mathbb{T} \ni t \mapsto n_{\mu,t}(M(t))$ is a complex $Q$-martingale and $n_{\mu,\bar{T}}(M(\bar{T})) = 1 + \int_0^{\bar{T}} n'_{\mu,t}(M(t)) \, dM(t)$. Since also $t \mapsto n'_{\mu,t}(M(t))$ is a $Q$-martingale it follows that

$$(A.60) \qquad n_{\mu,\bar{T}}(M(\bar{T})) = 1 + \int_0^{\bar{T}} E_Q(n'_{\mu,\bar{T}}(M(\bar{T})) | \mathcal{F}_t) \, dM(t).$$

Let $g \in C_0^\infty(\mathbb{R})$ be real-valued with Fourier transform $\hat{g}$. Multiplication of both sides of equality (A.60) with the complex number

$$c(\mu) = \frac{1}{\sqrt{2\pi}} e^{(-\mu^2/2)\langle\!\langle M,M\rangle\!\rangle(\bar{T})} \hat{g}(\mu)$$

gives

$$c(\mu) n_{\mu,\bar{T}}(M(\bar{T})) = c(\mu) + \int_0^{\bar{T}} E_Q(c(\mu) n'_{\mu,\bar{T}}(M(\bar{T}))|\mathcal{F}_t)\, dM(t).$$

Integration in $\mu$ and the stochastic Fubini theorem then give

$$(A.61) \qquad g(M(\bar{T})) = \int_{\mathbb{R}} c(\mu)\, d\mu + \int_0^{\bar{T}} E_Q(g'(M(\bar{T}))|\mathcal{F}_t)\, dM(t).$$

Since $(E_Q(g'(M(\bar{T}))|\mathcal{F}.))^2$ is a submartingale it follows that

$$E_Q\left(\int_0^{\bar{T}} (E_Q(g'(M(\bar{T}))|\mathcal{F}_t))^2\, d\langle\!\langle M,M\rangle\!\rangle(t)\right)$$
$$\leq E_Q((g'(M(\bar{T})))^2)\langle\!\langle M,M\rangle\!\rangle(\bar{T}),$$

which is finite. Therefore $\int_{\mathbb{R}} c(\mu)\, d\mu = E(g(M(\bar{T})))$, so (A.61) proves the representation formula (A.59) for $F \in C_0^\infty(\mathbb{R})$. The general case now follows by dominated convergence since $F$ in the lemma is the limit, in the topology defined by the norm $G \mapsto (E_Q(F(M(\bar{T}))^2 + F'(M(\bar{T}))^2))^{1/2}$, of a sequence in $C_0^\infty(\mathbb{R})$. $\qquad\square$

**Acknowledgments.** The authors would like to thank Nizar Touzi and Walter Schachermayer for fruitful discussions and for pointing out several references. We also thank the anonymous referees for constructive suggestions.

DEPARTMENT OF MATHEMATICS
UNIVERSITY OF BRITISH COLUMBIA
1984 MATHEMATICS ROAD
V6T 1Z2 CANADA
E-MAIL: ekeland@math.ubc.ca

EISTI
AVENUE DU PARC
95011 CERGY
FRANCE
E-MAIL: taflin@eisti.fr